\documentclass[11pt]{amsart}
\usepackage[utf8]{inputenc}
\usepackage[english]{babel}

\usepackage{amscd,amssymb}
\usepackage{amsmath}
\usepackage{amsthm}
\usepackage{graphicx}
\usepackage{subfigure,fancybox}
\usepackage{epic,eepic}
\usepackage{amstext}
\usepackage{esint}
\usepackage[square,numbers]{natbib}
\usepackage{booktabs}
\usepackage{hyperref}
\usepackage{comment}
\usepackage{mathtools}
\usepackage{tikz,pgfplots}
\usepackage{stmaryrd}

\def\elem{T}

\def\V{\mathcal{V}}

\def\Pnull{\mathbb{P}_0}

\def\T{\mathcal{T}}
\def\TH{\mathcal{T}_H}

\def\C{\mathbb{C}}

\def\PiH{\Pi_H}

\def\one{\mathbf{1}}

\def\nd{ \gamma_{\text{\scalebox{.65}[.8]{$\partial_n$}}}}

\newcommand{\tnorm}[2]{{\left\lVert #1 \right\rVert}_{L^2(#2)}}
\newcommand{\tnormf}[2]{{\| #1 \|}_{L^2(#2)}}

\newcommand{\vnormf}[2]{{\| #1 \|}_{\V_#2}}
\newcommand{\vnormof}[1]{{\| #1 \|}_{\V}}
\newcommand{\tsp}[3]{{\left( #1\,,\,#2 \right)}_{L^2(#3)}}
\newcommand{\tspf}[3]{{( #1\,,\,#2 )}_{L^2(#3)}}

\definecolor{myBlue}{RGB}{113,104,238} 
\definecolor{myGreen}{RGB}{114,175,30} 
\definecolor{myRed}{RGB}{180,50,50}  
\definecolor{myOrange}{RGB}{225,92,22}

\newtheorem{theorem}{Theorem}[section]

\newtheorem{lemma}[theorem]{Lemma}
\newtheorem{assumption}[theorem]{Assumption}

\newtheorem{conjecture}[theorem]{Conjecture}
\theoremstyle{definition}

\theoremstyle{remark}
\newtheorem{remark}[theorem]{Remark}
\numberwithin{theorem}{section}
\numberwithin{equation}{section}
\numberwithin{table}{section}
\numberwithin{figure}{section}

\allowdisplaybreaks

\usepackage[colorinlistoftodos,prependcaption]{todonotes}

\begin{document}
	
\begin{abstract}
	We propose a novel variant of the Localized Orthogonal Decomposition (LOD) method for time-harmonic  scattering problems of Helmholtz type with high wavenumber $\kappa$. On a coarse mesh of width $H$, the proposed method identifies local finite element source terms that yield rapidly decaying responses under the solution operator. They can be constructed to high accuracy from independent local snapshot solutions on patches of width $\ell H$ and are used as problem-adapted basis functions in the method. In contrast to the classical LOD and other state-of-the-art multi-scale methods, the localization error decays super-exponentially as the oversampling parameter $\ell$ is increased. This implies that optimal convergence is observed under the substantially relaxed oversampling condition $\ell \gtrsim (\log \tfrac{\kappa}{H})^{(d-1)/d}$ with $d$ denoting the spatial dimension. Numerical experiments demonstrate the significantly improved offline and online performance of the method also in the case of heterogeneous media and perfectly matched layers.
\end{abstract}
	
\title[]{Super-localized Orthogonal Decomposition for high-frequency Helmholtz problems}
\author[]{Philip Freese$^\dagger$, Moritz Hauck$^\dagger$, Daniel Peterseim$^\ddagger$}
\address{${}^{\dagger}$ Department of Mathematics, University of Augsburg, Universit\"atsstr.~12a, 86159 Augsburg, Germany}
\address{${}^{\ddagger}$ Department of Mathematics \& Centre for Advanced Analytics and Predictive Sciences (CAAPS), University of Augsburg, Universit\"atsstr.~12a, 86159 Augsburg, Germany}
\email{\{philip.freese, moritz.hauck, daniel.peterseim\}@math.uni-augsburg.de}
\thanks{The work of all authors is part of a project that has received funding from the European Research Council ERC under the European Union's Horizon 2020 research and innovation program (Grant agreement No.~865751).
The work of Philip Freese was partially funded by the Deutsche Forschungsgemeinschaft DFG (Project-ID 258734477 – SFB 1173)}
\maketitle

\vspace{1cm}

\section{Introduction}

 This paper studies the numerical solution of time-harmonic acoustic  scattering problems that can be modeled by the Helmholtz equation. The Helmholtz problem is an indefinite and non-hermitian problem and, especially for large wavenumbers $\kappa$, its numerical solution is a challenging task. Due to the highly oscillatory nature of the analytical solution and the $\kappa$-dependent pollution effect {\cite{BaS97}}, classical polynomial-based finite element methods need to meet very restrictive conditions on the mesh size $H$ of the underlying mesh. Typically, these conditions are much stronger than the minimal requirement $H \kappa  \lesssim  1$ from approximation theory needed for the approximation of an oscillatory function.

In the literature, there have been many attempts to tackle this issue. We highlight two classes of methods that are theoretically able to suppress the pollution effect, namely $hp$-finite elements and multi-scale methods. The strategy of $hp$-finite elements  \cite{Melenk2010hpfem,MeS11,Melenk2013hpfem} is to couple the polynomial degree $p$ of the approximation space to $\kappa$ in a logarithmic way. Then, the quasi-optimality of the numerical approximation can be ensured under the resolution condition $H \kappa  \lesssim p$. 

Multi-scale methods, in contrast, use problem-adapted ansatz spaces, which are constructed by solving multiple mutually independent local problems. An effective approach for constructing such problem-adapted ansatz spaces is the Localized Orthogonal Decomposition (LOD). The LOD was 
originally introduced as numerical homogenization method for elliptic diffusion problems with arbitrary rough coefficients \cite{MaP14,HeP13,MalP20} and later generalized to Helmholtz problems \cite{GaP15,Pet16,BGP17,Pet17,PeV20,HaPe21}. Its basis functions are computed by solving local subspace correction problems on element patches of size $\ell H$ with $\ell$ denoting the oversampling parameter of the method. Faithful numerical approximations to the Helmholtz problem are obtained under the resolution condition $H\kappa \lesssim 1$ and the oversampling condition $\ell \gtrsim \log \kappa$.  With their higher degree of adaptivity with respect to the problem, multi-scale methods and especially the LOD are intrinsically able to handle heterogeneous coefficients and singularities in the analytical solution. Those scenarios are not easily treated in  methods with universal shape functions such as $hp$-finite elements. 

Certainly, there exist many other numerical methods for the Helmholtz problem that were not mentioned yet. Widely known is the class of Trefftz methods \cite{trefftz0,trefftz3, trefftz1,trefftz}, which use, on each mesh element, functions that are locally solutions of the Helmholtz equation as test and trial functions (e.g. plain waves, generalized harmonic polynomials). Moreover, also other multi-scale methods like the multiscale finite element method (MsFEM) \cite{HoW97,EfH09} have been successfully applied to the Helmholtz problem \cite{Fu2017AFS,Fu2019WaveletbasedEM,Chen2021Helmh}.

The LOD implicitly computes its problem-adapted ansatz spaces as approximation to the space given by the application of the solution operator to some (coarse) classical finite element space; see \cite{Peterseim2021}. This connection has recently been exploited in \cite{HaPe21b} to develop a  conceptually new LOD multi-scale method for the elliptic multi-scale problem. It identifies local source terms in the coarse finite element space that yield rapidly decaying (in some cases even local) responses under the solution operator. This rapid decay makes it possible to approximate the global responses by localized counterparts, which are solutions to problems on element patches of size $\ell H$ in the coarse grid. These localized responses are then used as problem-adapted basis functions. The error caused by this approximation is henceforth referred to as localization error. For the elliptic multi-scale problem, the localization errors of the novel multi-scale method decays super-exponentially as $\ell$ is increased. This is a substantial improvement compared to the existing localization strategies  \cite{MaP14,HeP13,KorY16,KPY18,BrennerLOD} with exponentially decaying localization errors.

This paper aims to show that the novel localization strategy is not limited to the elliptic model problem but, just like the LOD, can be generalized to a large variety of problem classes beyond elliptic homogenization problems. As a proof of concept, this paper generalizes the novel localization strategy \cite{HaPe21b} to a class of indefinite and non-hermitian problems represented by the Helmholtz problem in the high-frequency regime. Under a stability assumption on the basis of the method and the resolution condition $H\kappa \ell \lesssim 1$ (in practice also $H\kappa \lesssim 1$ is sufficient), a $\kappa$-explicit stability and error analysis of the proposed multi-scale method is presented; see Theorems \ref{th:stab} and \ref{th:conv}. The stability and error estimates are explicit in quantities that are proportional to the smallest singular value of some patch-local coarse-scale operators. Numerical experiments clearly demonstrate a super-exponential decay of the singular values as $\ell$ is increased, although a rigorous mathematical proof is still open. Nevertheless, this motivates an ($\ell$-adaptive) a-posteriori error control strategy for the localization error using the easily accessible singular values as quasi-local error indicators. Given the super-exponential decay of the singular values, the oversampling condition needed for the stability of the method is $\ell \gtrsim (\log\tfrac{\kappa}{H})^{(d-1)/d}$ with $d$ denoting the spatial dimension. This is a major improvement compared to $\ell \gtrsim \log \kappa$ for the LOD. Under the asymptotically same condition, the proposed method yields an optimally convergent approximation which is, by the power $(d-1)/d$, better than the condition 
$\ell \gtrsim \log\tfrac{\kappa}{H}$ for the LOD.

This relaxed oversampling condition is of great practical importance as, especially for large $\kappa$, it enables significant computational savings compared to the LOD. First and foremost, the relaxed oversampling condition allows for smaller patch problems, which considerably reduces the offline computational costs. Moreover, also the online computational costs are lower due to the improved locality of the basis functions, which implies a sparser coarse system matrix.  

For the sake of simplicity, our analysis covers only Helmholtz problems in homogeneous media. Nevertheless, in the spirit of {\cite{BGP17,PeV20}, the method naturally extends} 
to the case of heterogeneous media which is demonstrated numerically in Section \ref{sec:numexp}.
In addition, this paper addresses the issue of physically meaningful boundary conditions. Although being widely used and convenient for mathematical theory, the impedance boundary condition as an approximation of the Dirichlet-to-Neumann map on some artificial boundary yields significant errors, especially for large $\kappa$; see \cite{Galkowski2021LocalAB}. Similarly, as in \cite{cgnt2018} for the LOD, we demonstrate that the proposed multi-scale method is naturally and easily combined with perfectly matched layers (PML) \cite{PML2,PML3,PML1} which are known to be an effective and efficient way to eliminate spurious reflections at the artificial boundary. 

The structure of this paper is as follows. Section \ref{sec:Model problem} briefly introduces the model problem and states some important analytical results. In Section \ref{sec:Prototypical multi-scale method}, we  present a prototypical multi-scale method  with optimal $\kappa$-independent convergence rates. Using the novel localization approach presented in Section \ref{sec:Novel localization strategy}, the method is then turned into a practically feasible method in Section \ref{sec:Novel super-localized multi-scale method}.
Finally, Section \ref{sec:numexp} illustrates the performance of the proposed method in numerical experiments. We show that the method can also be applied to heterogeneous media and is easily combined with the PML. 

\section{Model problem}\label{sec:Model problem}
Let us consider the Helmholtz equation with homogeneous impedance boundary conditions on a bounded polygonal Lipschitz domain $\Omega\subset \mathbb{R}^d$, $d = 1,2,3$ which is assumed to be scaled to unit size. Given a right-hand side $f\in L^2(\Omega)$, we seek $u\colon \Omega \rightarrow \mathbb C$ being the solution of
\begin{equation}\label{eq:classform}
\begin{aligned}
-\Delta u -\kappa^2u &= f\quad &&\text{ in }\Omega,\\
\nabla u \cdot n -i\kappa u &= 0 \quad &&\text{ on }\partial \Omega
\end{aligned}
\end{equation}
with $\kappa>0$ denoting the wavenumber, $i$ the imaginary unit, and $n$ the  outward unit normal vector. The weak formulation of \eqref{eq:classform} is based on the 
sesquilinear form  $a\colon \V\times\V\rightarrow \mathbb{C}$ which acts on the solution space $\V \coloneqq H^1(\Omega)$ (the space of complex-valued square-integrable functions on $\Omega$ with square-integrable weak derivative) and is defined as
\begin{equation*}
a(u,v) \coloneqq \tsp{\nabla u}{\nabla v}{\Omega} - \kappa^2\tsp{u}{v}{\Omega} - i\kappa\tsp{u}{v}{\partial \Omega}.
\end{equation*}
The inner products of the spaces $L^2(\Omega)$ and $L^2(\partial \Omega)$ (the spaces of complex-valued square-integrable functions on $\Omega$ and $\partial \Omega$, respectively) are denoted by
$\tsp{\cdot}{\cdot}{\Omega}$ and $\tsp{\cdot}{\cdot}{\partial \Omega}$, respectively. As usual in the Helmholtz context, the solution space $\V$ is endowed with the following $\kappa$-dependent norm
\begin{equation*}
\vnormof{u}^2 \coloneqq {\tnorm{\nabla u}{\Omega}^2 + \kappa^2 \tnorm{u}{\Omega}^2}.
\end{equation*}
With respect to this norm, the sesquilinear form $a$ is continuous, i.e., there is a $\kappa$-independent constant $C_a>0$ such that 
\begin{equation*}\label{eq:conta}
|a(u,v)|\leq C_a\, \vnormof{u}\vnormof{v}.
\end{equation*}
For a given $f\in L^2(\Omega)$, Fredholm theory \cite{Mel95} shows that there exists a unique weak solution $u \in \V$ to the weak formulation of \eqref{eq:classform}, such that, for all $v \in \V$, 
\begin{equation}\label{eq:wf}
a(u, v) = (f,v)_{L^2(\Omega)}
\end{equation}
satisfying 
\begin{equation}\label{polystab}
\vnormof{u} \leq C_{\mathrm{st}}(\kappa) \tnorm{f}{\Omega}
\end{equation}
with a $\kappa$-dependent constant $C_\mathrm{st}(\kappa)>0$. For the geometric configuration in this work, i.e., bounded Lipschitz domain and pure impedance boundary conditions, it can be shown that $C_\mathrm{st}$ depends polynomially on $\kappa$, i.e., $C_\mathrm{st} = \mathcal O(\kappa^n)$ for some $n\in \mathbb N$; see \cite{EsM12}. In the case where $\Omega$ is convex or smooth and star-shaped with respect to a ball, it is proved in \cite{Mel95} that $C_\mathrm{st} = \mathcal O(1)$. It should be noted that for more general boundary conditions and non-convex geometries, trapping scenarios can arise with $C_\mathrm{st}$ growing at least exponentially in $\kappa$; see  \cite{Betcke2011} for an example.

An immediate consequence of \eqref{polystab} is the inf-sup stability of $a$
\begin{equation}\label{eq:infsupcont}
0<\alpha(\kappa) \leq \adjustlimits{\inf}_{u\in \V} {\sup}_{v\in\V} \frac{\mathfrak{R} a(u,v)}{\vnormof{u}\vnormof{v}} = \adjustlimits{\inf}_{v\in \V} {\sup}_{u\in\V} \frac{\mathfrak{R} a(u,v)}{\vnormof{u}\vnormof{v}}
\end{equation}
with $\alpha(\kappa) = (2C_{\mathrm{st}}(\kappa)\kappa)^{-1}$ and $\mathfrak Rz$ denoting the real part of the complex number $z \in \mathbb C$.

We will refer to $\mathcal L \colon L^2(\Omega)\rightarrow \V$ as the solution operator of the Helmholtz problem, mapping $f\in L^2(\Omega)$ to the unique solution $u\in\V$ of the weak formulation \eqref{eq:wf}. For $z\in\C$ we denote with $\overline{z}$ the complex conjugate. The solution operator to the adjoint Helmholtz problem, i.e., $a(u,v)$ in \eqref{eq:wf} is substituted by $ \overline{a(v,u)}$, is denoted by $\mathcal L^*$. As shown in \cite[Lemma 3.1]{MeS11}, for $f\in L^2(\Omega)$ the solution operators $\mathcal L$ and $\mathcal L^*$ are connected by the relation
\begin{equation}\label{eq:relsolop}
\mathcal L^* f = \overline{\mathcal L \overline{\,f\,}}.
\end{equation}

\begin{remark}[Heterogeneous media, scatterers, boundary conditions]
	The construction of our method is easily extended to more complicated scenarios of heterogeneous materials, scattering problems and, for instance, mixtures of Dirichlet, Neumann and Robin boundary conditions. 
	It shall be mentioned, that the analytical well-posedness of Helmholtz problems in heterogeneous media is a delicate issue; see \cite{Sauter2018,Graham2018,Graham2019} for some recent results.
\end{remark}

\section{Prototypical multi-scale method}\label{sec:Prototypical multi-scale method}
This section introduces a prototypical multi-scale method which has convergence properties independent of the wavenumber. Its trial and test spaces are obtained by applying the (adjoint) solution operator to right-hand sides in classical finite element spaces. Such methods are only considered for theoretical purposes, as in general, global problems need to be solved for computing the method's basis functions. 

Let us introduce the possibly coarse, conforming mesh $\TH$ consisting of closed simplicial or quadrilateral elements with a diameter at most $H>0$. Henceforth, we assume that $\TH$ is  quasi-uniform and that its elements are non-degenerate in the sense of \cite[Definition 4.4.13]{BrS08}.
Denoting with $\mathbb P_0(\TH)$ the space of $\TH$-piecewise constant functions, we define  $\Pi_H\colon  L^2(\Omega)\rightarrow \Pnull(\TH)$ as the $L^2$-orthogonal projection onto $\Pnull(\TH)$. Recall that, for all $\elem\in \mathcal T_H$, it satisfies the following local stability and approximation properties (see  \cite{Poincareoriginal,Poincare})
\begin{equation}\label{e:L2proj}
\begin{aligned}
\|\PiH  v\|_{L^2(\elem)}&\leq \|v\|_{L^2(\elem)}\qquad &&\text{ for all }v\in L^2(\elem),\\
\|v-\PiH v\|_{L^2(\elem)}&\leq {\pi}^{-1}H \|\nabla v\|_{L^2(\elem)}\quad &&\text{ for all }v\in H^1(\elem).
\end{aligned}
\end{equation}

We introduce a prototypical problem-adapted Petrov--Galerkin method which is based on the trial and test spaces 
\begin{equation*}
\V_H \coloneqq \mathrm{span}\{\mathcal L \one_\elem\,|\,\elem\in \mathcal T_H\},\quad \V_H^* \coloneqq \mathrm{span}\{\mathcal L^*\one_\elem\,|\,\elem\in \mathcal T_H\},
\end{equation*}
where $\one_\elem$ denotes the characteristic function of an element $\elem\in \TH$. The prototypical multi-scale method  seeks $u_H\in \V_H$ such that, for all $v_H \in \V_H^*$,
\begin{equation}\label{e:galerkinideal}
a(u_H,v_H) = \tspf{f}{v_H}{\Omega}.
\end{equation}

The following necessary assumption ensures that the coarse mesh $\TH$ is able to resolve oscillatory Helmholtz solutions; see Lemma \ref{l:ua} below.
\begin{assumption}[Resolution condition]\label{ass:mesh}
	Suppose that the mesh size of $\TH$ satisfies
	\begin{equation*}\label{eq:minH}
	H\kappa \leq \frac{\pi}{\sqrt{2}}.
	\end{equation*}
\end{assumption}
Under this assumption, the well-posedness and the $\kappa$-independent approximation properties of the prototypical Petrov--Galerkin method can be proved.
\begin{lemma}[Stability and $\kappa$-independent approximation]\label{l:ua}
	If Assumption \ref{ass:mesh} is  fulfilled, then the sesquilinear form $a$ is inf-sup stable with regard to the trial space $\V_H$ and the test space $\V_H^*$, i.e., there exists $C_\mathrm{id}>0$ independent of $\kappa$ and $H$ such that
	$$C_\mathrm{id} \adjustlimits \inf_{u_H\in \V_H}\sup_{v_H \in \V_H^*}\frac{\mathfrak R a(u_H,v_H)}{\vnormof{u_H}\vnormof{v_H}} \geq  \alpha(\kappa) >0$$
	Here, $\alpha$ denotes the inf-sup constant of the continuous problem \eqref{eq:wf}. 
	
	Moreover, there exists $C_\mathrm{er}>0$ independent of $\kappa$ and $H$ such that, for all right-hand sides $f \in H^s(\Omega)$ with $s \in [0,1]$, the unique solution $u_H$ of the Petrov--Galerkin method \eqref{e:galerkinideal} satisfies the $\kappa$-independent error bound 
	\begin{equation*}
	\frac{\pi}{2}\,{\vnormof{u - u_H}} \leq H \tnormf{f -\PiH f}{\Omega}\leq  C_\mathrm{er} H^{1+s}{\|f\|_{H^s(\Omega)}}.
	\end{equation*}
\end{lemma}
\begin{proof}
	For a proof see, e.g., \cite[Theorem 3.9 and Example 3.10a]{Peterseim2021} or \cite[Lemma 4.6 and 4.7]{HaPe21}. 
\end{proof}

\section{Localization strategy}\label{sec:Novel localization strategy}
The canonical basis functions $\{\mathcal L\one_\elem\,|\,\elem \in \TH\}$ and $\{\mathcal L^* \one_\elem\,|\, \elem \in \TH\}$ of the problem-adapted trial and test spaces $\V_H$ and  $\V_H^*$, respectively, are non-local and have a slow (algebraic) decay. For a practically feasible variant of method \eqref{e:galerkinideal}, localized bases have to be identified. Recently, in \cite{HaPe21b}, a novel localization approach was introduced for the elliptic model problem which has superior localization properties compared to other state-of-the-art approaches  \cite{MaP14,HeP13,KorY16,KPY18,BrennerLOD}. This section extends this novel localization approach to a class of indefinite non-hermitian problems, using the Helmholtz problem as an example. 
In contrast to \cite{HaPe21}, the target regime of the proposed method is the high-frequency case where the oscillatory behavior of the solution is just resolved by the coarse mesh $\TH$. 
For this regime, the multi-resolution approach in \cite{HaPe21} is not applicable.

By relation \eqref{eq:relsolop}, it suffices to analyze the localization of the trial space $\V_H$; a localized basis of the test space can then be obtained without further computation. The idea of the localization strategy is to identify local $\TH$-piecewise constant source terms that yield rapidly decaying (or even local) responses under the solution operator $\mathcal L$ of the Helmholtz problem. 

\begin{figure}
	\includegraphics[width=.275\linewidth]{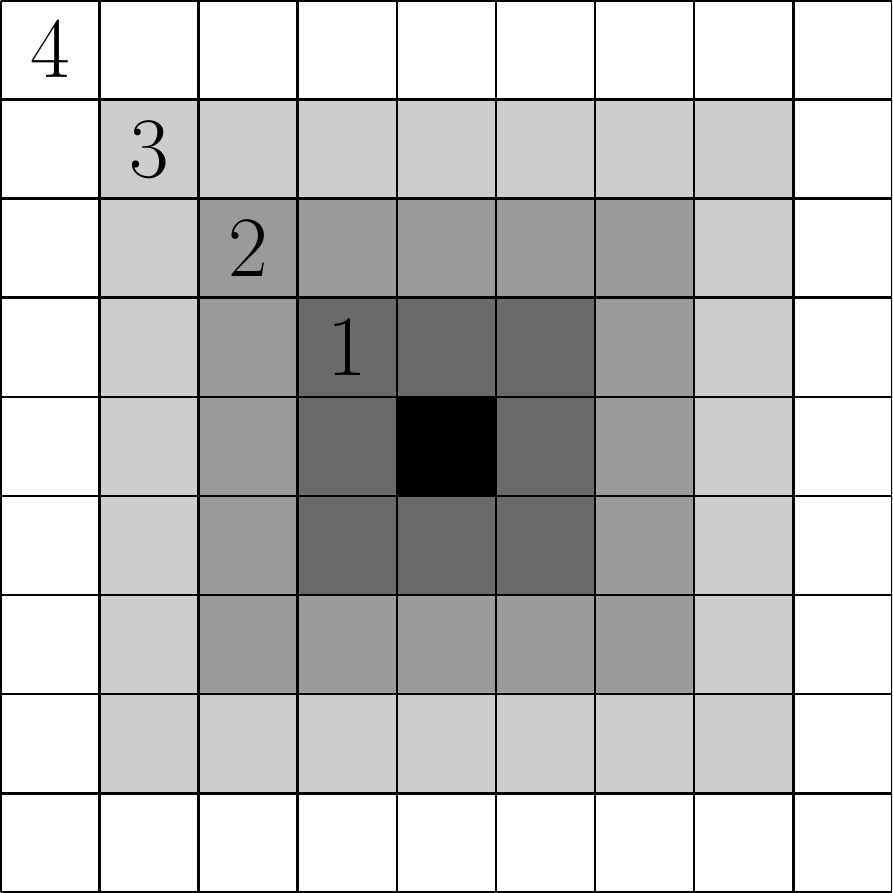}
	\caption{Illustration of a $\ell$-th order patch for $\ell = 1,\dots,4$ with gray scale indicating  the order. }
	\label{fig:patch}
\end{figure}

Our localization is based on (local) patches, given as neighborhoods of mesh elements, in the coarse mesh $\TH$. The first order element patch $\mathsf{N}(S)=\mathsf{N}^1(S)$ of a union of elements $S\subset \Omega$ is given by
\begin{equation*}
	\mathsf{N}^1(S)\coloneqq \bigcup \left\{\elem\in \TH\,|\,\elem\cap S\neq \emptyset\right\}.
\end{equation*}
The $\ell$-th order patch $\mathsf{N}^\ell(T)$, $\ell=2,3,4,\dots$, of $T$ is then recursively given by
\begin{equation*}
\mathsf{N}^\ell(T)\coloneqq \mathsf{N}^{1}(\mathsf{N}^{\ell-1}(T));
\end{equation*}
see Figure \ref{fig:patch} for a schematic illustration.

For the subsequent derivation, we fix an oversampling parameter $\ell \in \mathbb N$ and denote, for arbitrary $T \in \TH$, the  $\ell$-th order patch of $T$ by $\omega\coloneqq \mathsf{N}^\ell(\elem)$. Let  $\ell$ be chosen such that no patch coincides with the domain $\Omega$. On the patch $\omega$, we define the space  $\V_\omega\coloneqq \{v|_\omega\,|\,v \in \V\}$, i.e.,  the restriction of $\V$ to the patch $\omega$. Furthermore, let $\T_{H,\omega}$  denote the submesh of $\TH$ with elements in $\omega$ and let $\Pi_{H,\omega}\colon  L^2(\omega) \rightarrow \Pnull(\T_{H,\omega})$ denote  the $L^2(\omega)$-orthogonal projection onto $\Pnull(\T_{H,\omega})$. 

The novel localization strategy requires the solution of local patch problems with Dirichlet and impedance boundary conditions prescribed on $\Gamma \coloneqq\partial \omega \backslash \partial \Omega$ and $\partial \omega \cap \partial \Omega$, respectively. For proving the well-posedness of such problems, we need the Friedrich's inequality (see e.g. \cite{BrS08,R2013}) which states the existence of a constant $C_\mathrm{F}>0$ such that, for all $v \in \V_{\omega,\Gamma} \coloneqq \{v \in \V_\omega\,|\, v|_\Gamma = 0 \}$, 
\begin{equation}\label{eq:Fr}
\tnormf{v}{\omega} \leq C_\mathrm{F}H\ell \tnormf{\nabla v}{\omega}.
\end{equation}
Henceforth, we suppose that $C_\mathrm{F}$ is independent of the parameters $H$ and $\ell$; see also the following remark.

\begin{remark}[Friedrich's constant]
	Provided that, for example, $\TH$ is a Cartesian mesh, all patches are intervals/rectangles/cuboids of diameter $\mathcal O(H\ell)$. By transformation techniques, one can then prove that $C_\mathrm{F}$ is independent of $H$ and $\ell$. For the general setting, however, the dependence of $C_\mathrm{F}$ on the parameter $\ell$, cannot be figured out explicitly. Nevertheless, if the meshes $\TH$ are generated by uniform refinement  of an initial coarse mesh $\T_0$ (capturing the geometric characteristics of the domain $\Omega$), it is reasonable to assume that $C_\mathrm{F}$ is well-behaved. This can be justified by the fact that the geometric complexity of the mesh and its element patches is already determined by the initial mesh and does not increase by uniform refinement.
	
\end{remark}
The following assumption poses a stronger condition on the smallness of the mesh size than Assumption \ref{ass:mesh}.
\begin{assumption}[Resolution condition  revisited]\label{ass:meshsize}
	Suppose that the mesh size of $\TH$ satisfies
	\begin{equation*}
	H\kappa\ell \leq \frac{1}{C_\mathrm{F}\sqrt{2}}.
	\end{equation*}
\end{assumption}

In practice, the additional $\ell$-dependence seems artificial, as for the LOD \cite{Pet17,HaPe21}, such an additional condition is not needed. Indeed, the numerical experiments in Section \ref{sec:numexp} indicate that, in practice, the weaker Assumption \ref{ass:mesh} is sufficient.

The following lemma states the coercivity of the  sesquilinear form $a_\omega\colon \V_\omega \times \V_\omega \rightarrow \mathbb C$ 
\begin{equation*}
a_\omega(u,v) \coloneqq \tspf{\nabla u}{\nabla v}{\omega} - \kappa^2 \tspf{u}{v}{\omega} - i\kappa \tspf{u}{v}{\partial \Omega\cap \partial \omega}
\end{equation*}
with respect to the norm  $\vnormf{\cdot}{\omega}^2 \coloneqq {\tnormf{\nabla \cdot}{\omega}^2 + \kappa^2 \tnormf{\cdot}{\omega}^2}$.
\begin{lemma}[Coercivity of $a_\omega$]\label{lemma:coerc}
	If Assumption \ref{ass:meshsize} is fulfilled, it holds, for all $v \in \V_{\omega,\Gamma}$, that
	\begin{equation*}
	\mathfrak Ra_\omega(v,v) \geq \frac{1}{3}\vnormf{ v}{\omega}^2.
	\end{equation*}
\begin{proof}
	For all $v \in \V_{\omega,\Gamma}$, we obtain using  Friedrich's inequality \eqref{eq:Fr} and Assumption \ref{ass:meshsize}
	\begin{equation*}
	\mathfrak Ra_\omega(v,v) 
	= \tnormf{\nabla v}{\omega}^2 - \kappa^2 \tnormf{v}{\omega}^2
	\geq (1-C_\mathrm{F}^2H^2\ell^2\kappa^2)\tnormf{\nabla v}{\omega}^2 \geq \frac{1}{2}\tnormf{\nabla v}{\omega}^2.
	\end{equation*}
	The coercivity follows, utilizing that, for all $v \in \V_{\omega,\Gamma}$,
	\begin{equation*}
	\vnormf{v}{\omega}^2 =  \tnormf{\nabla v}{\omega}^2 +\kappa^2  \tnormf{v}{\omega}^2 \leq (1+C_\mathrm{F}^2\kappa^2H^2\ell^2)\tnormf{\nabla v}{\omega}^2 \leq \frac{3}{2}\tnormf{\nabla v}{\omega}^2.\qedhere
	\end{equation*} 
\end{proof}
\end{lemma}

We aim to identify  (almost) local basis functions $\varphi = \varphi_{\elem,\ell}\in \V_H$ associated with the element $\elem \in \TH$. For coefficients $(c_K)_{K \in \T_{H,\omega}}$ that need to be determined afterwards, the construction of these basis functions follows the ansatz
\begin{equation*}
	\varphi = \mathcal L g \quad \text{ with  }\quad g = g_{\elem,\ell} \coloneqq \sum_{K\in \mathcal T_{H,\omega}}c_K \one_K.
\end{equation*}
The Galerkin projection of such a given $\varphi$ onto the local subspace $\V_{\omega,\Gamma}$ defines a localized approximation $\varphi^\mathrm{loc} = \varphi^\mathrm{loc}_{\elem,\ell} \in \V_{\omega,\Gamma}$, which satisfies, for all $v\in \V_{\omega,\Gamma}$,
\begin{equation}\label{e:patchproblem}
a_\omega(\varphi^\mathrm{loc},v)= \tspf{g}{v}{\omega}.
\end{equation}
Due to the coercivity of $a_\omega$ (Lemma \ref{lemma:coerc}) and the Lax--Milgram lemma, $\varphi^{\mathrm{loc}}$ is well-defined. In general, the local function $\varphi^{\mathrm{loc}}$ is a poor approximation of the possibly global function $\varphi$. However, there exist nontrivial choices of $g$ that yield highly accurate approximations in the $\V$-norm. The following discussion requires a brief reminder on traces of $\V_\omega$-functions; see \cite{LiM72a} for details. We denote the trace operator restricted to $\Gamma\subset \partial \omega$ by
\begin{equation*}
\gamma_{0} = \gamma_{0,\omega}\colon \V_\omega \rightarrow H^{1/2}(\Gamma)
\end{equation*}
and abbreviate its range by $X\coloneqq H^{1/2}(\Gamma)$. For given $t \in X$, we define the Helmholtz-harmonic extension $\gamma_{0}^{-1}$ which is a continuous right-inverse of $\gamma_{0}$. To be precise, the extension $\gamma_{0}^{-1}t$ satisfies $\gamma_0 \gamma_{0}^{-1}t= t$ and is Helmholtz-harmonic, i.e., for all $v \in \V_{\omega,\Gamma}$,
\begin{equation}\label{eq:pp}
	a_\omega(v,\gamma_{0}^{-1}t)= 0.
\end{equation}
Note that well-posedness is again ensured by the coercivity of $a_\omega$; see Lemma \ref{lemma:coerc}. The normal derivative of $\varphi^\mathrm{loc}$, denoted by $\nd \varphi^{\mathrm{loc}}$, is defined as the element of $X^\prime$ that satisfies, for all $v \in\V_\omega$, 
\begin{equation*}
	\langle \nd \varphi^\mathrm{loc}\,,\,\gamma_0v\rangle_{X^\prime\times X} = -\tspf{g}{v}{\omega} +a_\omega(\varphi^\mathrm{loc},v).
\end{equation*}
Now we state our main observation, which establishes a criterion for the choice of $g$ ensuring the smallness of the localization error. Using the properties of the trace and Helmholtz-harmonic extension operators, we obtain, for all $v \in \V$,
\begin{equation}\label{e:nd1}
	a(\varphi-\varphi^\mathrm{loc} ,v) = (g,v)_{L^2(\omega)}-a_\omega(\varphi^\mathrm{loc} ,v)=  -\langle \nd \varphi^\mathrm{loc}, \gamma_{0}\, v\rangle_{X^\prime\times X}= (g, \gamma_{0}^{-1} \gamma_{0}\, v)_{L^2(\omega)},
\end{equation}
where we used the Helmholtz-harmonic extension in the last equality. Hence, by the inf-sup stability of the continuous problem \eqref{eq:infsupcont}, a small localization error is equivalent to a small norm of the normal derivative $\nd \varphi^\mathrm{loc}$  which, in turn, is equivalent to choosing  $g$ (almost) $L^2$-orthogonal to the space
 \begin{equation}\label{e:Y}
	Y\coloneqq \gamma_{0}^{-1}X\subset \V_\omega
\end{equation}
of Helmholtz-harmonic functions on $\omega$. An optimal realization of $g$ can hence be achieved by the singular value decomposition of the operator $\Pi_{H,\omega}|_Y$, i.e., choosing $g$ as the right singular vector corresponding to the smallest singular value.

Henceforth, we suppose that the above choice of $g$ satisfies, for some  parameter \linebreak  $\sigma_\elem(\kappa,H,\ell)>0$, the following estimate
\begin{equation}\label{e:orth}
\|g\|_{Y^\prime}\coloneqq 
\sup_{v \in Y} \frac{\tspf{g}{v}{\omega}}{\vnormf{v}{\omega}} \leq \sigma_\elem(\kappa,H,\ell)\|g\|_{\mathcal V_{\omega}^\prime}.
\end{equation}

The quantity $\sigma_T$ coincides, up to a constant, with the smallest singular value of the operator $\Pi_{H,\omega}|_{Y}$ and will be used in the remainder as a measure for the (quasi-) orthogonality of $g$ on $Y$. The dependence of $\sigma_T$ on the wavenumber is due to the wavenumber-dependent space $Y$ and the norm of the solution space.

We conjecture that $\sigma_T$ decays super-exponentially in $\ell$ which is subsequently justified with a numerical experiment. 

\begin{conjecture}[Super-exponential decay]\label{con:decay}
	The quantity $\sigma_T$ decays super-exponentially in $\ell$, i.e., there exist constants $C_\mathrm{sd}(\kappa,H,\ell)>0$  depending polynomially on $\kappa, H$, and $\ell$, but being independent of $T$ and $C>0$ independent of $\kappa,H,\ell,$ and $T$ such that
	\begin{equation*}
	\sigma_{\elem}(\kappa,H,\ell) \leq C_\mathrm{sd}(\kappa,H,\ell)  \exp\left(-C\ell^{\frac{d}{d-1}}\right).
	\end{equation*}
\end{conjecture}

\begin{figure}[h]
	\includegraphics[width=.47\textwidth]{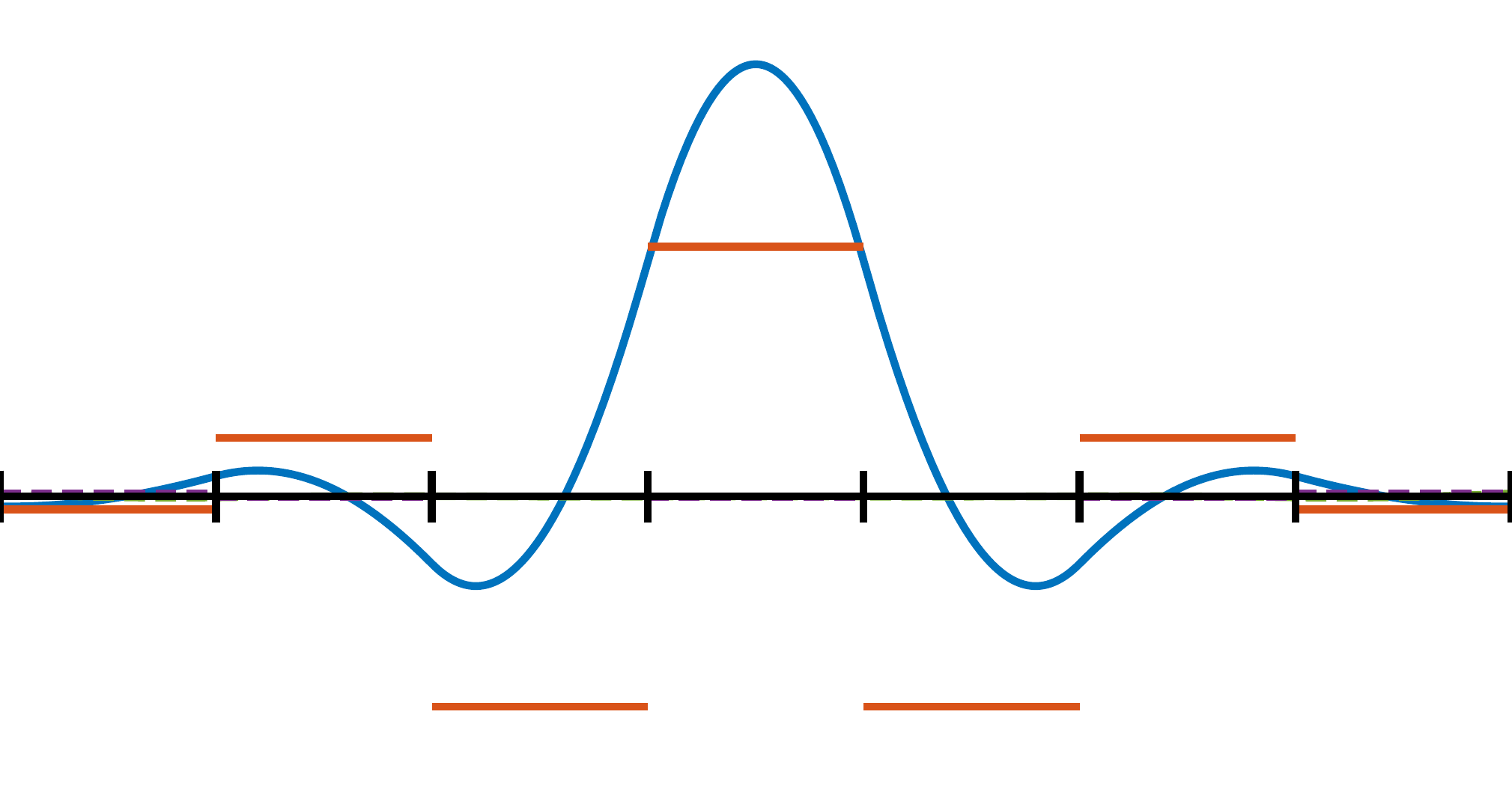}\hfill
	\includegraphics[width=.47\textwidth]{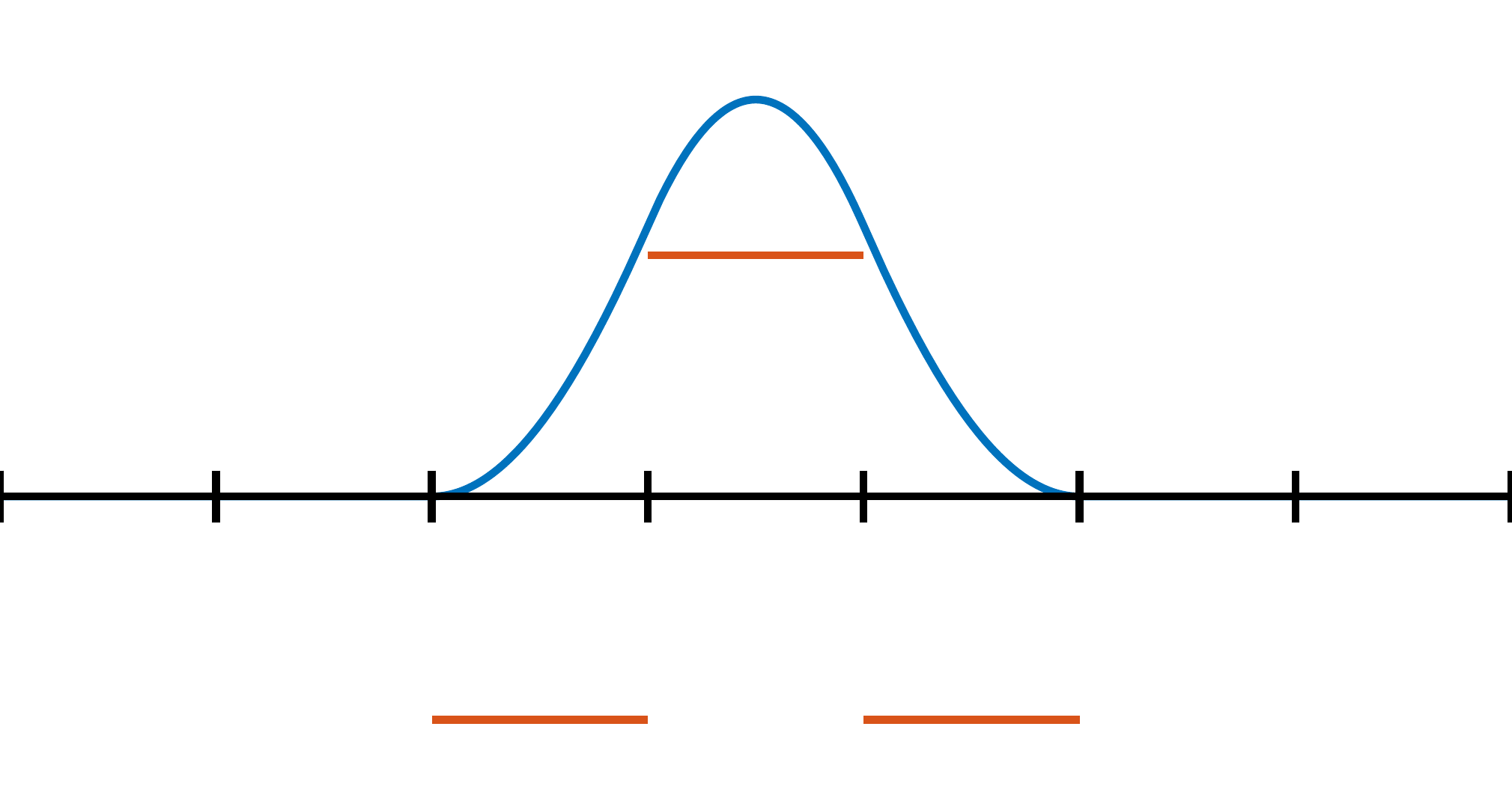}
	\includegraphics[width=.47\textwidth]{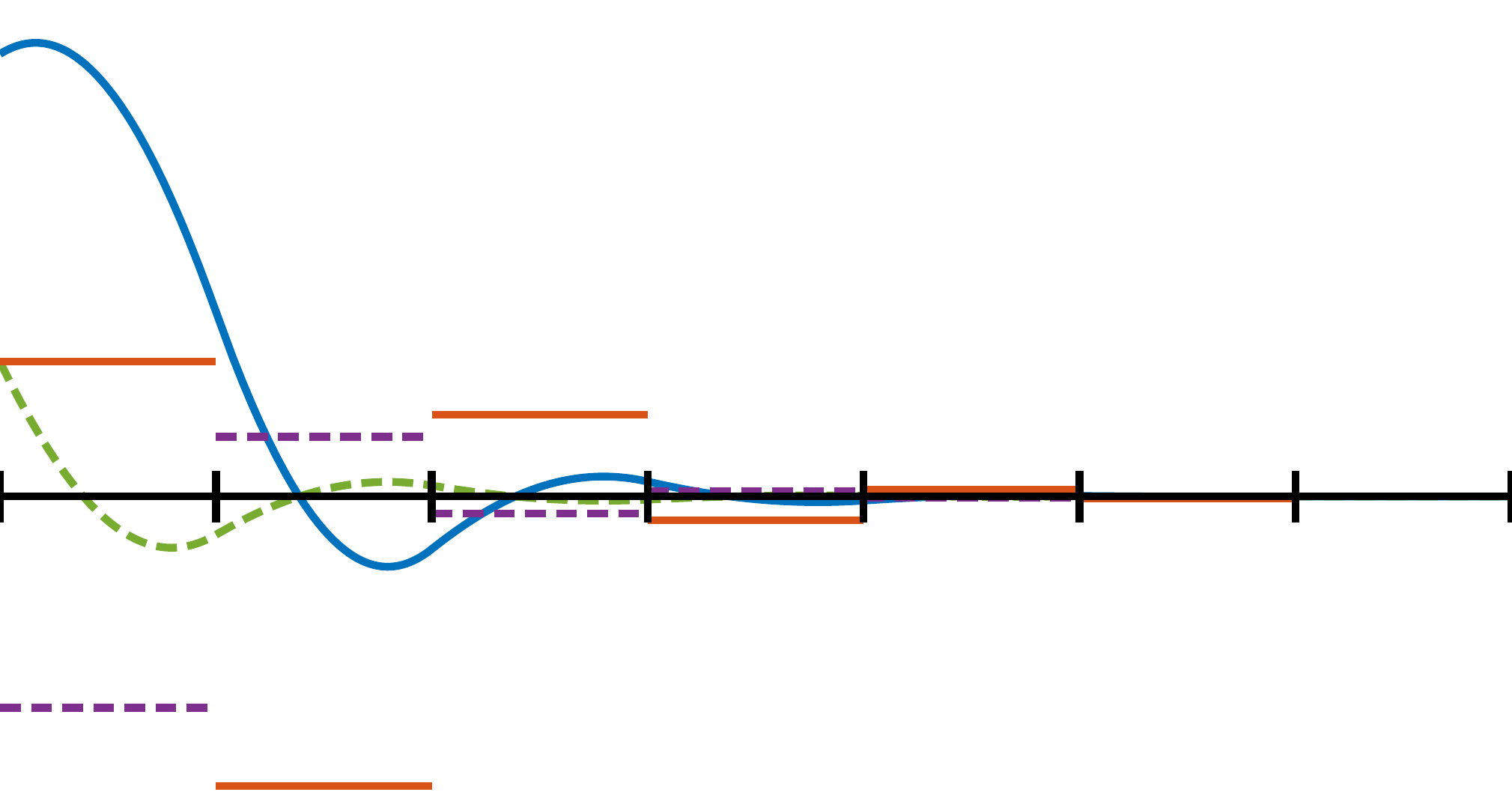}\hfill
	\includegraphics[width=.47\textwidth]{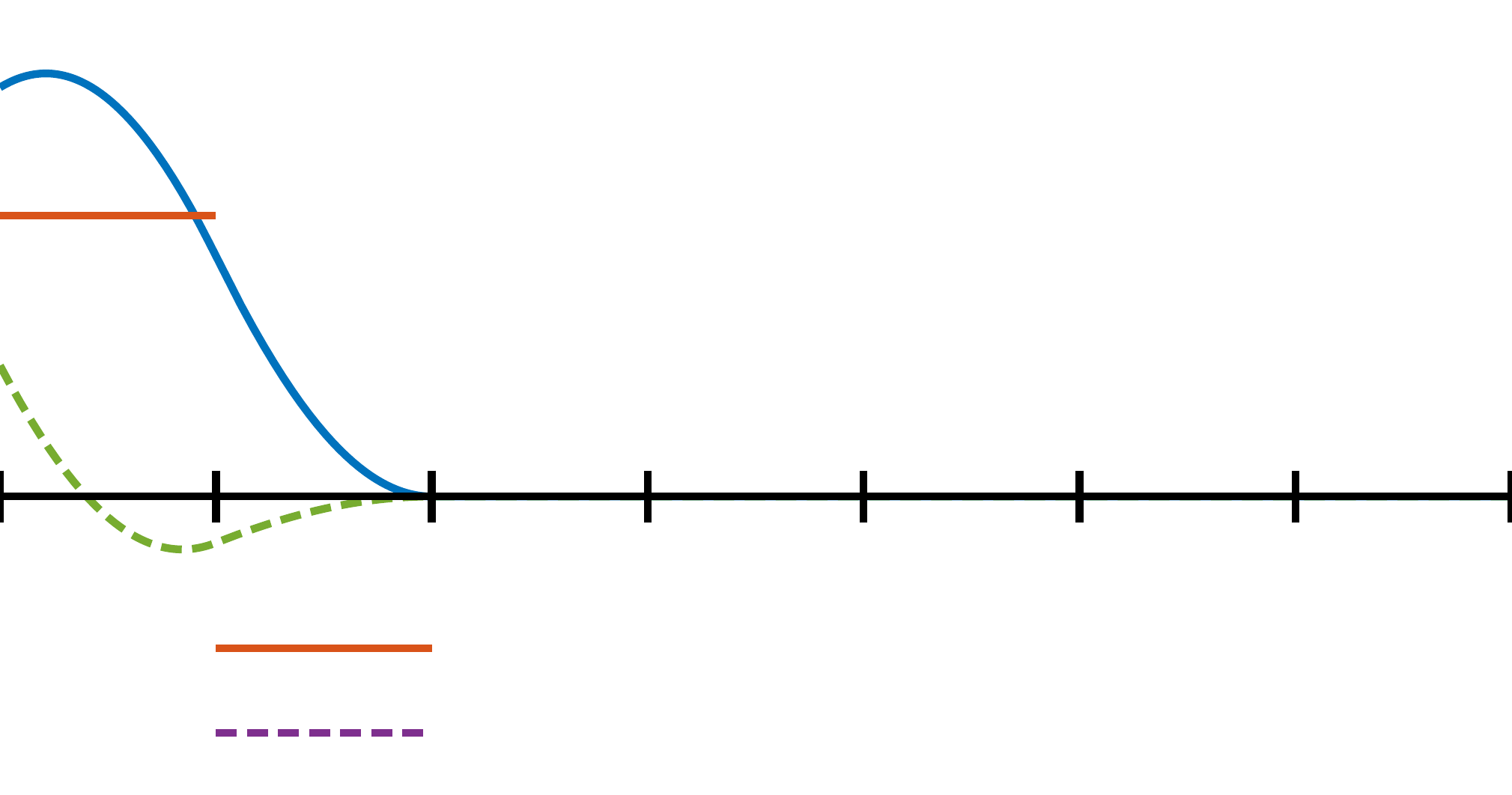}
	\caption{Global ideal LOD basis (left) and local SLOD basis for $\ell = 1$ (right) with their corresponding $L^2$-normalized right-hand sides $g$ in one space dimension for an interior element (top) and an element at the boundary (bottom). The real (resp. imaginary) parts are depicted using solid (resp. dashed) lines. }
	\label{fig:basis1d}
\end{figure}

\begin{remark}[The case $d=1$]
	For one spatial dimension, the space $Y$ of Helmholtz-harmonic functions is at most two-dimensional. Thus, for $\ell\geq 1$, $g$ can indeed be chosen  $L^2$-orthogonal on $Y$, i.e., the basis is in fact local. This locality is in line with our conjecture interpreting $\tfrac{d}{d-1}$ as infinity. Figure \ref{fig:basis1d} shows the local SLOD basis functions and the ideal (non-localized) LOD basis functions from \cite{HaPe21} corresponding to an interior element and an element at the boundary. It shall be noted that the imaginary part of the LOD basis function corresponding to the interior element is very small and thus not visible. For the respective SLOD basis function, the imaginary part is zero.
\end{remark}

\begin{figure}[h]
	\includegraphics[width=.49\textwidth]{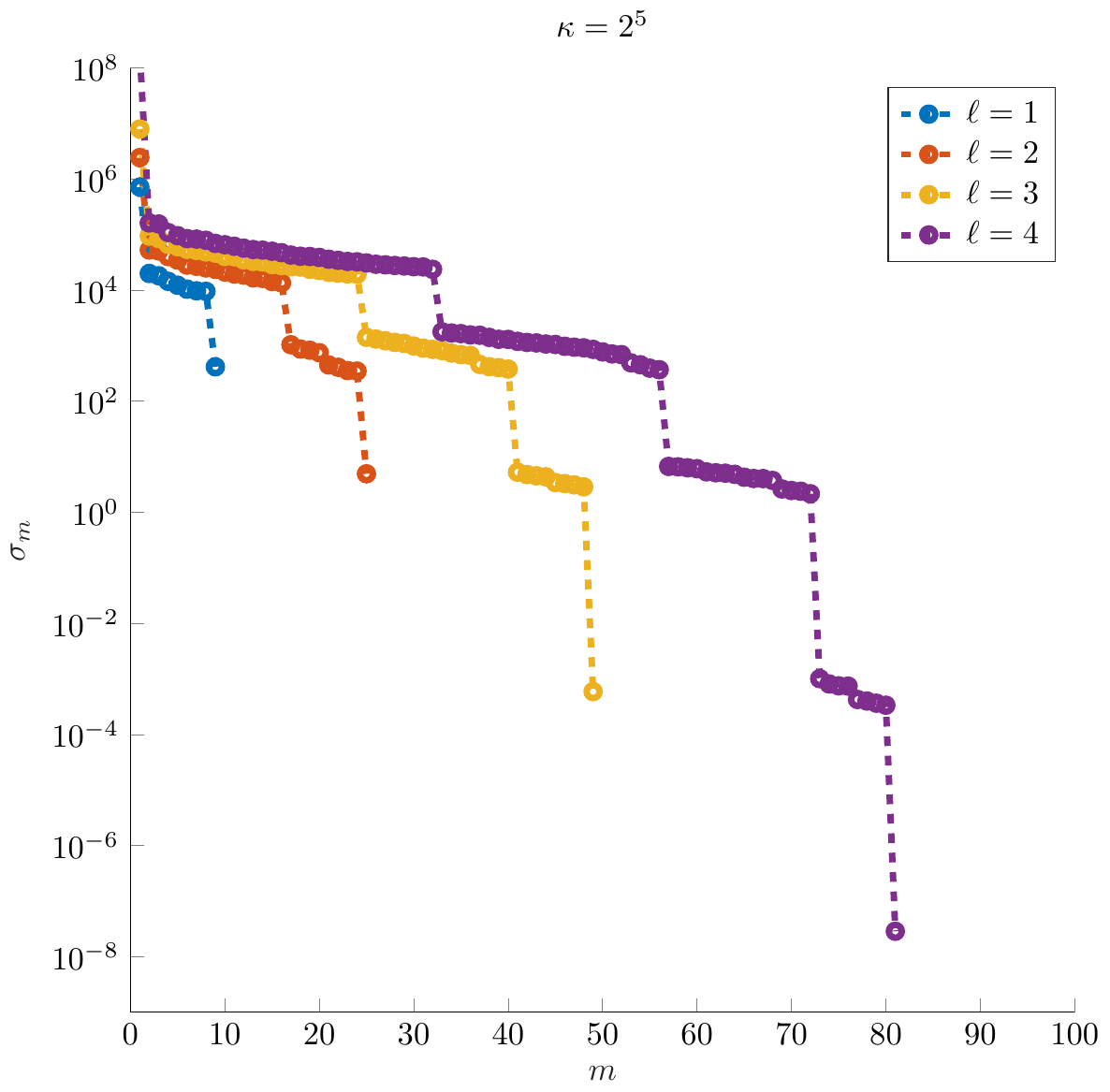}
	\includegraphics[width=.49\textwidth]{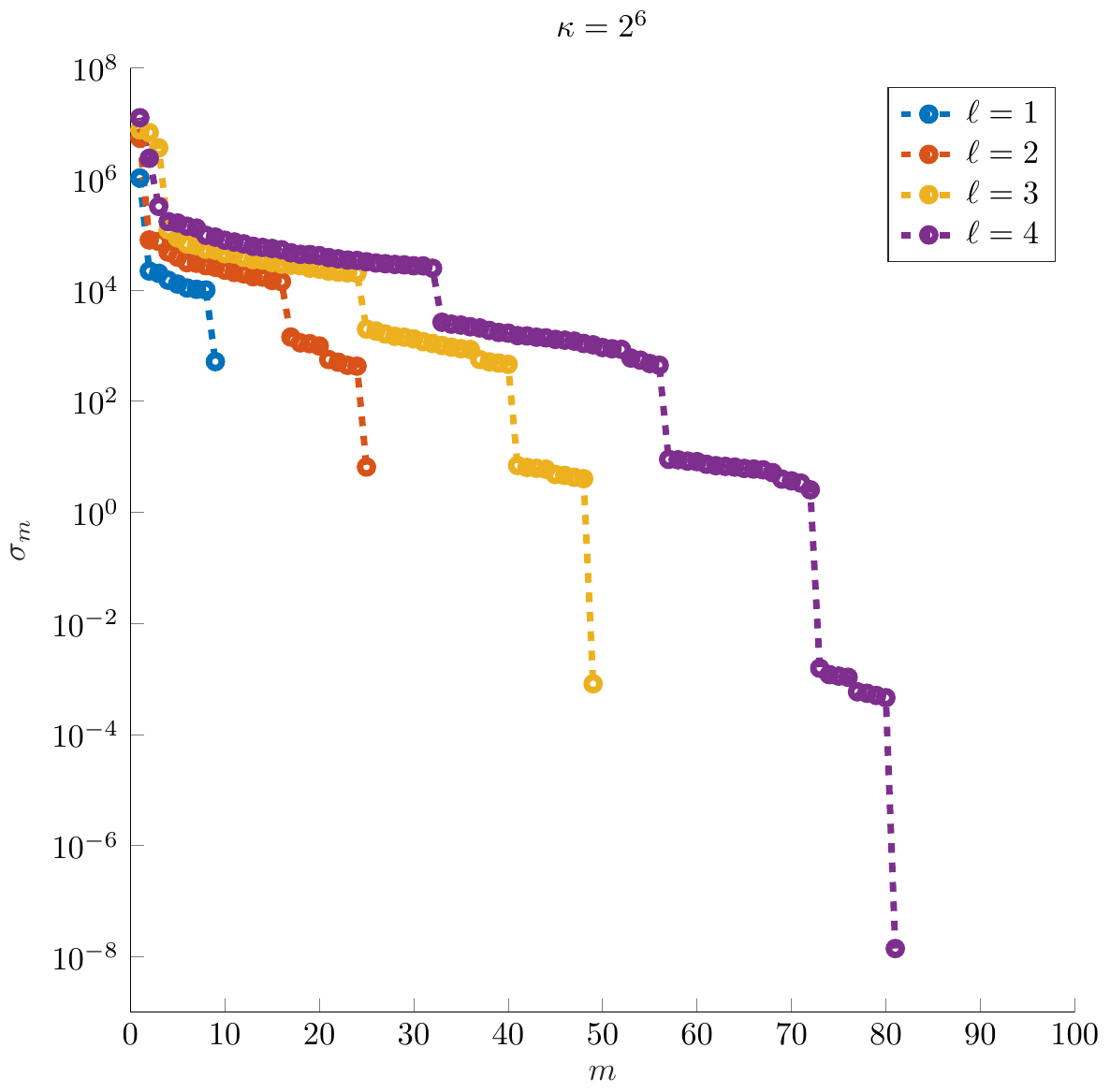}
	\caption{Singular values $\sigma_m$ of the operator $\Pi_{H,\omega}|_{Y}$ in descending order for several oversampling parameters for $\kappa = 2^5$ (left) and $\kappa = 2^6$ (right). }
	\label{fig:singular_values}
\end{figure}

Next, we provide numerical experiments that demonstrate the super-exponential decay of the singular values of the operator $\Pi_{H,\omega}|_{Y}$ and thus justify  Conjecture \ref{con:decay}. We select an element $T$ of a fixed Cartesian mesh $\TH$ which is far away from the boundary and consider the patches $\omega = \mathsf N^\ell(T)$, $\ell = 1,\dots,4$; see Section \ref{sec:numexp} for the precise setup of the numerical experiment. Figure \ref{fig:singular_values} depicts the singular values of the operator $\Pi_{H,\omega}|_{Y}$  in a semi-logarithmic plot.

One observes that, for all $\ell$, the decay of the singular values is step-like with jumps between the plateaus doubling from left to right. Note that, for $\ell = 4$, the last jump is smaller than expected which is due to the finite precision arithmetic of the computer. This observation is precisely the super-exponential decay (quadratic exponentially for $d = 2$) of the singular values. Let us point out that, for all $\ell$, the smallest singular value is well separated. Hence, it is clear that $g$, cf.\,\eqref{e:orth}, should be chosen as the right singular vector corresponding to this singular value. For patches touching the boundary, one observes a similar decay behavior as in Figure \ref{fig:singular_values}.

Using techniques from LOD theory \cite{MalP20,Peterseim2021}, one can rigorously derive a pessimistic bound for  $\sigma_T$. This pessimistic bound  shows that the novel localization techniques performs at least as good as state-of-the-art approaches like the LOD.

\begin{lemma}[Pessimistic exponential decay]
	The quantity $\sigma_T(\kappa,H,\ell)$ decays at least exponentially in $\ell$, i.e., there exist constants $C, C^\prime >0$ independent of $\kappa,H,\ell,$ and $T$ such that 
	\begin{equation*}
	\sigma_T(\kappa,H,\ell) \leq C' \max\{1,H^{-1}\kappa^{-1}\} \exp\left(-C\ell\right).
	\end{equation*}
\end{lemma}
\begin{proof}
	This result can be proved using techniques from LOD theory and can be obtained by straight-forward modifications of the proof of  \cite[Lemma 6.4]{HaPe21b}.
\end{proof}
 
\section{Super-localized multi-scale method}\label{sec:Novel super-localized multi-scale method}

Using the novel localization strategy of the previous Section, we turn the prototypical multi-scale method \eqref{e:galerkinideal} into a feasible scheme. For a fixed oversampling parameter $\ell$, we define the ansatz space of the localized method as the span of the localized basis functions $\varphi^\mathrm{loc}_{\elem,\ell}$ defined in \eqref{e:patchproblem}, i.e., 
\begin{equation*}
\V_{H,\ell} \coloneqq \mathrm{span}\{\varphi_{\elem,\ell}^\mathrm{loc}\,|\,\elem\in \mathcal T_H\},\quad \V_{H,\ell}^* \coloneqq \mathrm{span}\{\varphi_{\elem,\ell}^{\mathrm{loc},*}\,|\,\elem\in \mathcal T_H\}.
\end{equation*}
The localized method then  calculates the Petrov--Galerkin approximation with trial space $\V_{H,\ell}$ and 
 test space $\V_{H,\ell}^*$, i.e.,  it 
	 seeks $u_{H,\ell}\in \V_{H,\ell}$ such that, for all $v_{H,\ell}\in \V_{H,\ell}^*$,
\begin{equation}\label{e:wflocmethod}
a(u_{H,\ell},v_{H,\ell}) = \tspf{f}{v_{H,\ell}}{\Omega}.
\end{equation}

A minimal requirement for the stability and convergence of the Galerkin method \eqref{e:wflocmethod} is that $\{g_{\elem,\ell}\,|\,\elem \in \TH\}$ spans $\Pnull (\mathcal T_H)$ in a stable way. Numerically, this can be ensured as outlined in \cite[Appendix B]{HaPe21b}. For the subsequent numerical analysis, we make the following assumption.
\begin{assumption}[Riesz stability]\label{a:Rieszbasis} 
	The set $\{g_{\elem,\ell}\,|\,\elem \in \TH\}$ is a Riesz basis of  $\Pnull(\mathcal T_H)$, i.e., there exists a constant  $C_{\mathrm{rb}}(\kappa,H,\ell)>0$ depending only polynomially on $H, \ell$ such that, for all $(c_\elem)_{T \in \TH}$, it holds
	\begin{equation*}
	C^{-1}_{\mathrm{rb}}(\kappa,H,\ell)\sum_{\elem \in \mathcal T_H} |c_\elem|^2 \|g_{T,\ell}\|_{\mathcal V_{\omega}^\prime }^2 \leq \Big\lVert \sum_{\elem \in \mathcal T_H} c_\elem g_{\elem,\ell}\Bigl\rVert_{\mathcal V^\prime }^2 \;\leq C_{\mathrm{rb}}(\kappa,H,\ell) \sum_{\elem \in \mathcal T_H} |c_\elem|^2\|g_{T,\ell}\|_{\mathcal V_\omega^\prime}^2,
	\end{equation*}
	where $\omega = \mathsf N^\ell(T)$.
\end{assumption}

In what follows, we investigate the inf-sup stability and the convergence properties of the novel super-localized multi-scale method. The respective estimates are explicit in the quantity
\begin{equation*}
\sigma(\kappa,H,\ell) \coloneqq \max_{\elem\in\TH}\sigma_T(\kappa,H,\ell).
\end{equation*}

The following theorem shows that the novel method is inf-sup stable under a condition on the oversampling parameter $\ell$. 
\begin{theorem}[Stability]\label{th:stab}
	Let  Assumption \ref{ass:meshsize} and \ref{a:Rieszbasis} be satisfied and let $\ell$ be chosen such that
	\begin{equation}\label{eq:condinfsup}
	\varepsilon(\kappa,H,\ell)\coloneqq
	\alpha^{-1}(\kappa)(1+3C_a)C_a C_\mathrm{ol}\,\ell^{d/2} C_\mathrm{rb}^{1/2}(\kappa,H,\ell)\sigma(\kappa,H,\ell)\leq \min\left\{\frac{1}{2},\frac{\alpha(\kappa)}{27C_aC_\mathrm{id}}\right\}
	\end{equation}
	 with $\alpha$ denoting the inf-sup constant of the continuous problem \eqref{eq:infsupcont} and $C_\mathrm{ol}>0$ depending only on mesh properties of $\TH$ (non-degenerateness and quasi-uniformity).  Then, the localized method \eqref{e:wflocmethod} is inf-sup stable, i.e., there exists a constant $C_\mathrm{lo}>0$ independent of $\kappa$, $H$, and $\ell$ such that
	  \begin{equation*}
	  C_\mathrm{lo} \adjustlimits \inf_{u_{H,\ell} \in \V_{H,\ell}} \sup_{v_{H,\ell} \in \V_{H,\ell}^*} \frac{\mathfrak R a(u_{H,\ell},v_{H,\ell})}{\vnormof{u_{H,\ell}}\vnormof{v_{H,\ell}}} \geq  \alpha(\kappa).
	  \end{equation*}
\end{theorem}
\begin{proof}
	Let us consider the following bijective mapping between $\V_H$ and $\V_{H,\ell}$
	\begin{equation*}
	\iota\colon 
	\V_H \rightarrow \V_{H,\ell},\; u_H \coloneqq  \sum_{\elem \in \mathcal T_H}c_T \varphi_{\elem,\ell} \mapsto \sum_{\elem \in \mathcal T_H}c_T \varphi_{\elem,\ell}^\mathrm{loc} \eqqcolon u_{H,\ell}.
	\end{equation*}
	\emph{First}, we show the continuity of $\iota$. For arbitrary $u_{H,\ell} \in \V_ {H,\ell}$, the triangle inequality yields  
	\begin{align*}\label{eq:triu}
	\vnormof{u_{H,\ell}}\leq \vnormof{u_{H}} + \vnormof{u_H-u_{H,\ell}}.
	\end{align*}
	Using the inf-sup stability of the continuous problem \eqref{eq:infsupcont} combined with \eqref{e:nd1}, we obtain for the second term
	\begin{align*}
	\alpha(\kappa) \vnormof{u_H-u_{H,\ell}} &\leq \sup_{\vnormof{v} = 1} {\mathfrak R a(u_H-u_{H,\ell},v)}
	\leq \sup_{\vnormof{v} = 1}\big\vert \sum_{\elem \in \mathcal T_H}c_T \tspf{g_{T,\ell}}{\gamma_{0}^{-1}\gamma_{0}v}{\omega}\big\vert\\
	&\leq \sup_{\vnormof{v}=1 }\sum_{\elem \in \mathcal T_H}\sigma_T(\kappa,H,\ell) |c_\elem|\,\|g_{T,\ell}\|_{\V^\prime_\omega}\vnormf{\gamma_{0}^{-1}\gamma_{0}v}{\omega},
	\end{align*}
	where we used \eqref{e:orth} in the last inequality. Lemma \ref{lemma:coerc} yields the estimate  $\vnormf{\gamma_{0}^{-1}\gamma_0 v}{\omega} \leq (1+3C_a)\vnormf{v}{\omega}$  independent of the patch $\omega$. The finite overlap  of the patches yields
	\begin{align*}
		\sum_{\elem \in \mathcal T_H} \vnormf{v}{\omega}^2 \leq \sum_{\elem \in \mathcal T_H} C_{\mathrm{ol}}^2 \ell^d \vnormf{v}{\elem}^2 = C_{\mathrm{ol}}^2 \ell^d \vnormof{v}^2,
	\end{align*}
	where $C_\mathrm{ol}^2\ell^d$ bounds the number of patches containing a fixed mesh element. 	From this and the Cauchy--Schwarz inequality, we get
	\begin{align*}
	\alpha(\kappa) \vnormof{u_H - u_{H,\ell}}
	&\leq (1+3C_a)C_\mathrm{ol} \ell^{d/2}\sigma(\kappa,H,\ell)\sqrt{\sum_{\elem \in \mathcal T_H}|c_T|^2\|g_{T,\ell}\|^2_{\V^\prime_\omega}} \\
	&\leq (1+3C_a) C_\mathrm{ol}\ell^{d/2} C_\mathrm{rb}^{1/2}(\kappa,H,\ell) \sigma(\kappa,H,\ell) \Big\|\sum_{\elem \in \mathcal T_H}c_Tg_{T,\ell}\Big\|_{\V^\prime}\\
	&\leq (1+3C_a)C_a C_\mathrm{ol}\ell^{d/2} C_\mathrm{rb}^{1/2}(\kappa,H,\ell)\sigma(\kappa,H,\ell) \vnormof{u_H}
	\end{align*}
	and thus, using \eqref{eq:condinfsup}, $\vnormof{u_{H,\ell}}\leq \tfrac{3}{2}\vnormof{u_{H}}$, i.e., the continuity of $\iota$.
	Similarly, one can show $\vnormof{u_{H}}\leq 2\vnormof{u_{H,\ell}}$, i.e., the continuity of $\iota^{-1}$. The same  estimates can analogously be shown for 
	\begin{equation*}
	\iota^*\colon 
	\V_H^* \rightarrow \V_{H,\ell}^*,\; v_H \coloneqq  \sum_{\elem \in \mathcal T_H}c_T \varphi_{\elem,\ell}^* \mapsto \sum_{\elem \in \mathcal T_H}c_T \varphi_{\elem,\ell}^{\mathrm{loc},*} \eqqcolon v_{H,\ell}.
	\end{equation*}
	
	\emph{Second}, we show the inf-sup stability of the localized problem using the inf-sup stability of the continuous problem \eqref{eq:infsupcont}. We consider a  fixed but arbitrary $u_{H,\ell} \in \V_{H,\ell}$ and define $u_H \coloneqq \iota^{-1}u_{H,\ell} \in \V_H$. Furthermore, we set $v_{H,\ell} \coloneqq  \iota^{*} v_H  \in \V_{H,\ell}^*$, where $v_H \in \V_H^*$ is chosen such that
	\begin{equation*}
	   {\mathfrak R a(u_{H},v_H)}\geq \frac{\alpha(\kappa)}{C_\mathrm{id}}{\vnormof{u_{H}}\vnormof{v_H}},
	\end{equation*}
	cf.  Lemma \ref{l:ua}. Algebraic manipulations and elementary estimates yield 
	\begin{align*}
	\mathfrak R a(u_{H,\ell},v_{H,\ell}) \geq  \mathfrak R a(u_H,v_H) -|a(u_{H,\ell}-u_H,v_H)|-|a(u_{H,\ell},v_{H,\ell}-v_H)|.
	\end{align*}
	We estimate the terms on the right-hand side separately. For the first term, we obtain using the continuity of $\iota$ and $\iota^*$ that
	\begin{equation*}
	\mathfrak R a(u_H,v_H) \geq \frac{ \alpha(\kappa)}{C_\mathrm{id}} \vnormof{u_H}\vnormof{v_H}\geq \frac{\alpha(\kappa)}{3C_\mathrm{id}}\vnormof{u_{H,\ell}}\vnormof{v_{H,\ell}}.
	\end{equation*}
	For the second term, one obtains using \eqref{eq:condinfsup} 
	\begin{equation*}
	|a(u_{H,\ell}-u_H,v_H)|\leq C_a \vnormof{u_{H,\ell}-u_H}\vnormof{v_H} \leq \frac{\alpha(\kappa)}{9C_\mathrm{id}} \vnormof{u_{H,\ell}}\vnormof{v_{H,\ell}}.
	\end{equation*}
	The third term can be estimated analogously. Absorbing the second and the third term in the first, the inf-sup stability of the localized method follows.
\end{proof}

Utilizing Conjecture \ref{con:decay}, one can rewrite the oversampling condition \ref{eq:condinfsup} as  
\begin{equation}\label{eq:relos}
\ell \gtrsim (\log \tfrac{\kappa}{H})^{(d-1)/d}.
\end{equation}
As seen in the following theorem, the same  asymptotically condition also guarantees an optimal order of convergence. In both cases, this is a substantial improvement compared to the oversampling conditions for the LOD, which are $\ell \gtrsim \log \kappa$ for stability and $\ell \gtrsim \log \tfrac{\kappa}{H}$ for an optimal order of convergence. It shall be noted that the $\kappa$-dependence in the above oversampling conditions can be eliminated using Assumption \ref{ass:mesh} or \ref{ass:meshsize}.

\begin{theorem}[Convergence]\label{th:conv}
	Let the assumptions from Theorem \ref{th:stab} be fulfilled. Then, the solution $u_{H,\ell}$ of the localized Petrov--Galerkin approximation \eqref{e:wflocmethod} satisfies
	\begin{equation}\label{eq:conv}
	\vnormof{u-u_{H,\ell}} \leq \frac{2}{\pi} C_\mathrm{er} H^{1+s} \|f\|_{H^s(\Omega)} + \delta(\kappa,H,\ell)\tnormf{f}{\Omega}
	\end{equation}
	with $$\delta(\kappa,H,\ell) = 3\big(1+C_aC_\mathrm{id}\alpha^{-1}(\kappa)\big)C_\mathrm{id}\kappa^{-1}\alpha^{-1}(\kappa)\varepsilon(\kappa,H,\ell)$$
	and $C_\mathrm{er}$ from Lemma \ref{l:ua}.
\end{theorem}
\begin{proof}
	For this proof, we use the notation from the proof of Theorem \ref{th:stab}. We begin estimating with the triangle inequality
	\begin{equation*}
	\vnormof{u-u_{H,\ell}} \leq \vnormof{u-u_H} + \vnormof{u_H- u_{H,\ell}}.
	\end{equation*}
	The first term can be estimated using Lemma \ref{l:ua}. For the second term, we apply Strang's lemma \cite[Lemma 2.25]{ErG04} as $u_{H,\ell} \in \V_{H,\ell}$ can be seen as non-conforming and non-consistent approximation to $u_H \in \V_H$
	\begin{align*}
		\vnormof{u_H-u_{H,\ell}} &\leq \big(1+C_a C_\mathrm{id} \alpha^{-1}(\kappa) \big)\inf_{w_{H,\ell}\in \V_{H,\ell}} \vnormof{u_H-w_{H,\ell}} \\&\qquad +C_\mathrm{id} \alpha^{-1}(\kappa) \sup_{v_{H,\ell} \in \V_{H,\ell}^*}\frac{|a(u_H,v_{H,\ell})-\tspf{f}{v_{H,\ell}}{\Omega}|}{\vnormof{v_{H,\ell}}}.
	\end{align*}
	Here, the first term can be estimated, choosing $w_{H,\ell} \coloneqq \iota u_H$. Using similar arguments as in the proof of Theorem \ref{th:stab}, this yields
	\begin{align*}
	 \vnormof{u_H-\iota u_{H}} \leq \varepsilon(\kappa,H,\ell)\vnormof{u_H}\leq C_\mathrm{id}\kappa^{-1}\alpha^{-1}(\kappa)\varepsilon(\kappa,H,\ell)\tnorm{f}{\Omega}.
	\end{align*}
	For the second term, elementary algebraic manipulations yield, for all $v_H \in \V_H^*$,
	\begin{equation*}
	a(u_H,v_{H,\ell})-\tspf{f}{v_{H,\ell}}{\Omega} = \tspf{f}{v_H-v_{H,\ell}}{\Omega} - a(u_H,v_H-v_{H,\ell}).
	\end{equation*}
	Choosing $v_H \coloneqq \iota^{*,-1}v_{H,\ell}$, we obtain
	\begin{align*}
	&|a(u_H,v_{H,\ell})-\tspf{f}{v_{H,\ell}}{\Omega}|\\&\quad \leq 2\kappa^{-1}\varepsilon(\kappa,H,\ell)\tnormf{f}{\Omega}\vnormof{v_{H,\ell}} + 2C_a C_\mathrm{id}\kappa^{-1}\alpha^{-1}(\kappa)\varepsilon(\kappa,H,\ell)\tnormf{f}{\Omega}\vnormof{v_{H,\ell}}.
	\end{align*}
	Putting together the estimates finishes the proof.
\end{proof}

\begin{remark}[A-posteriori error control strategy for localization error]\label{r:posteriori}
Whilst the first term on the right-hand side of \eqref{eq:conv} can be controlled a-priori, an a-posteriori error control strategy for the second term (the localization error) seems beneficial.
For all $T \in \TH$, the singular values of the patch-local operators $\Pi_{H,\omega}|_Y$  might be used as (quasi-)local error indicators. It shall be noted that no extra computations are needed for evaluating these error indicators. In case that the error indicator shows a large local error, one might locally increase the oversampling parameter $\ell$ and thereby decrease the localization error. This approach yields a numerical algorithm guaranteeing a prescribed upper bound for the localization error.
\end{remark}

\section{Numerical experiments}\label{sec:numexp}
In this section, we investigate the proposed multi-scale method (henceforth referred to as SLOD) numerically and compare it with the LOD from \cite{HaPe21}. All our experiments were done using Matlab and a short version of the code is available as a supplement. It shall be noted that this code is for demonstration purposes only. Hence, it is not optimized and, e.g., does not exploit the structure of the coefficient.   

In the subsequent numerical experiments, we consider the domain $\Omega \coloneqq (0,1)^2$ endowed with a coarse Cartesian mesh $\TH$ of mesh size $H$. 
For discretizing the continuous patch-problems \eqref{e:patchproblem} and \eqref{eq:pp}, we use the $\mathcal Q_1$-finite element method on fine Cartesian meshes of the respective patches with mesh size~$2^{-10}$. All errors are computed against a $\mathcal Q_1$-finite element reference solution on the global  Cartesian mesh of mesh size $2^{-10}$.\\
The space $Y$ from  \eqref{e:Y} of Helmholtz harmonic functions on the patch $\omega = \mathsf N^\ell (T)$ is sampled using $5\cdot \# \T_{H,\omega}$ samples of random discrete Dirichlet data on $\partial \omega \backslash\partial \Omega$. The random Dirichlet data is generated by linearly interpolating independently and identically distributed (iid) values (from a uniform distribution) prescribed at the boundary vertices of the fine Cartesian patch-mesh lying on $\partial \omega \backslash\partial \Omega$. For more information on efficient sampling techniques for spaces of harmonic functions, see also \cite{BuS18}.
Next, the computed reduced space is used for calculating the singular value decomposition of the operator $\Pi_{H,\omega}|_{Y}$ in order to determine $g_{\elem,\ell}$ as well as its corresponding singular value $\sigma_{\elem}$ which might be used as error indicator for the a-posteriori error control strategy; see Remark \ref{r:posteriori}. It shall be noted that, for patches $\omega$ that are close to the boundary $\partial\Omega$, the choice of $g$ is more involved. Stability for such patches can be ensured by allowing the communication between at most $\mathcal O(\ell^d)$ patches. The corresponding algorithm can be found in the supplementary material; for a discussion of the algorithm, see \cite[Appendix B]{HaPe21b}.  

\begin{remark}[Computational costs]
The computation can be divided into an offline phase being independent of the right-hand side and an online phase which needs to be repeated if the right-hand side is changed. In the offline phase, the basis functions of the method are precomputed and the coarse stiffness matrix is assembled. Note that (for homogeneous media) it suffices to compute only $\mathcal O(\ell^d)$ basis functions, whilst the remaining ones can be obtained by translation; see \cite[Section 3]{GaP15}.
For a fixed oversampling parameter $\ell$, the computational costs for the proposed multi-scale method are comparable to those of the LOD (considering the usual LOD implementation described in \cite{EHM19}). However, due to the relaxed oversampling condition \eqref{eq:relos}, significantly smaller $\ell$'s are sufficient for reaching a prescribed level of accuracy. This shrinks the computational costs, as both, the overall number of patch-problems and their size is reduced considerably. 
For the online phase, the smaller oversampling parameter allows for a sparser system matrix which, in turn, reduces the computational costs.
\end{remark}

\subsection{Super-exponential decay of localization error}\label{sec:Super-exponential decay of localization error}
\begin{figure}[h]
	\includegraphics[width=.49\textwidth]{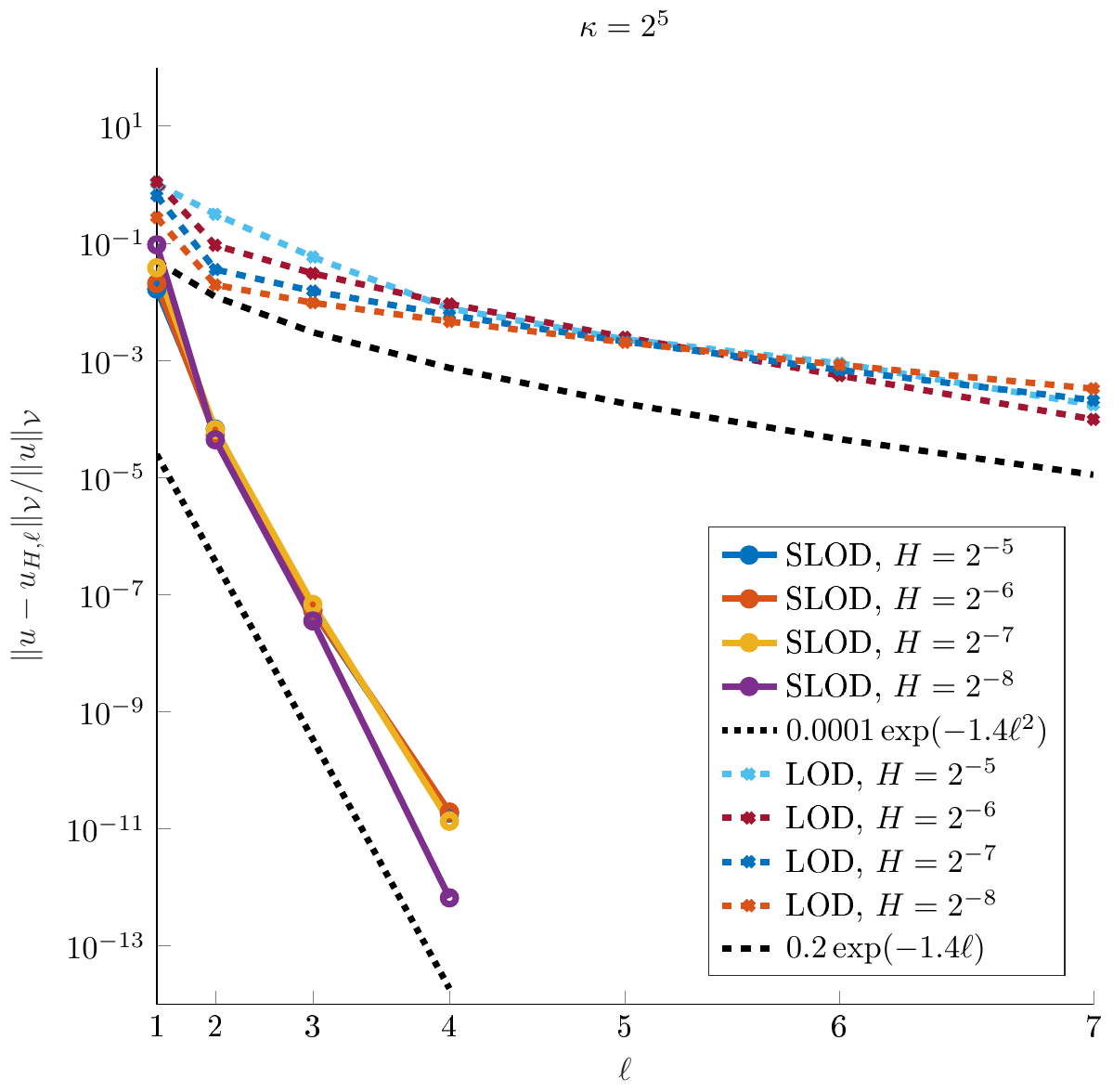}
	\includegraphics[width=.49\textwidth]{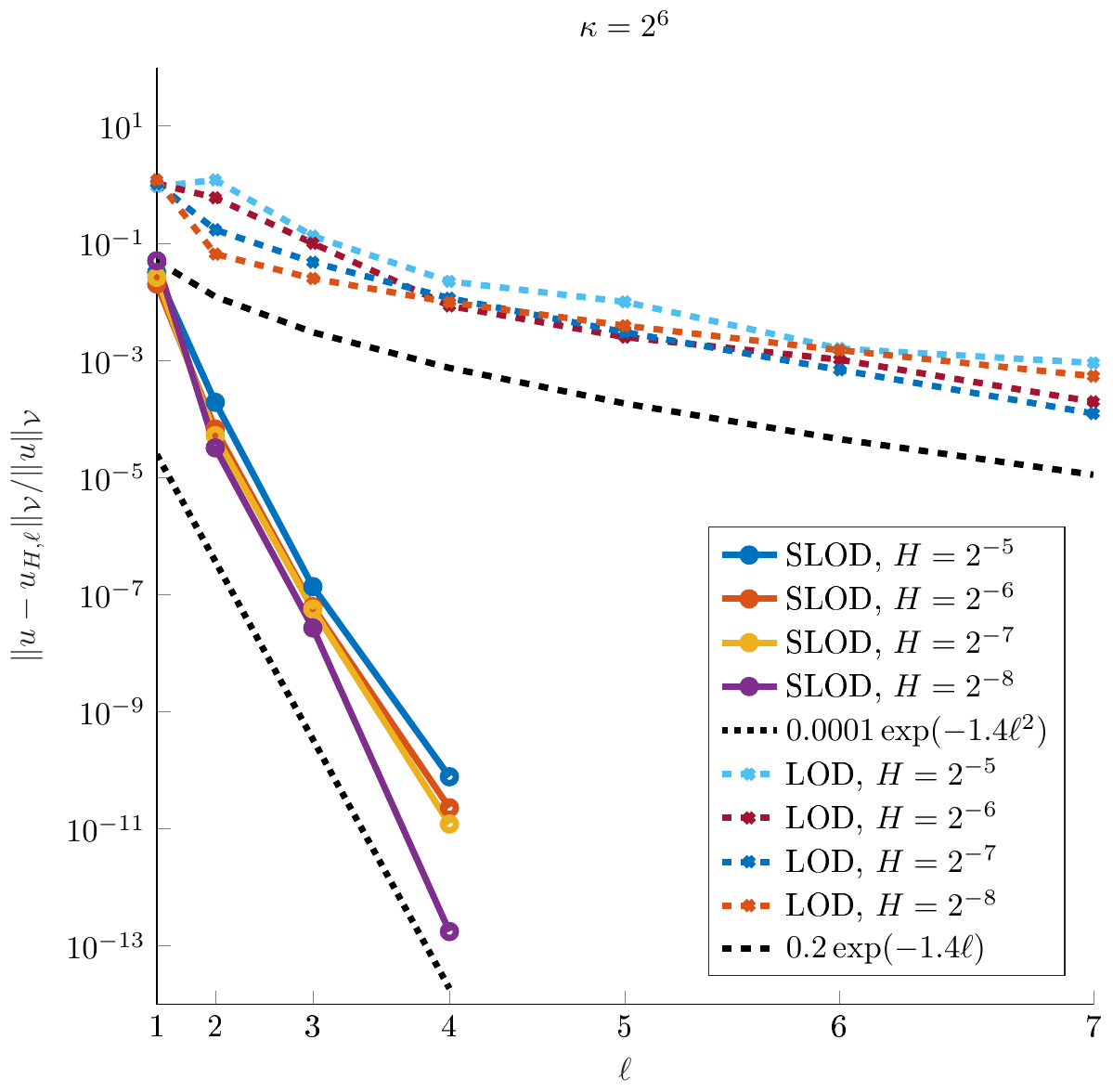}
	\caption{Localization errors for the SLOD and LOD for several values of $H$ for $\kappa = 2^5$ (left) and $\kappa = 2^6$ (right). }
	\label{fig:loc}
\end{figure}

For this numerical experiment, we consider the right-hand side $f \equiv 1$, as for this choice the first term in \eqref{eq:conv} vanishes and thus, only the localization error $\delta$ remains. 
Figure \ref{fig:loc} shows the relative $\V$-norm localization errors of our localization approach (referred to as
SLOD) and for the stabilized LOD from \cite{HaPe21}. The localization errors are plotted for several coarse grids $\TH$ in dependence of $\ell$ for $\kappa =  2^5, 2^6$. As reference, we indicate lines showing the expected rates of decay of the localization errors.  Please note the special scaling of the axes which is chosen such that quadratically exponentially decaying functions appear linear with negative slope. 

Figure \ref{fig:loc} numerically confirms the super-exponentially decay of the localization errors of the proposed multi-scale method. The localization error of the LOD decays exponentially; see \cite{HeP13,MaP14}. This numerical experiment confirms, that the SLOD has localization errors several orders of magnitude smaller than  the LOD.

\subsection{Optimal convergence under mesh refinement}\label{sec:Optimal convergence under mesh refinement}
\begin{figure}[h]
	\includegraphics[width=.49\textwidth]{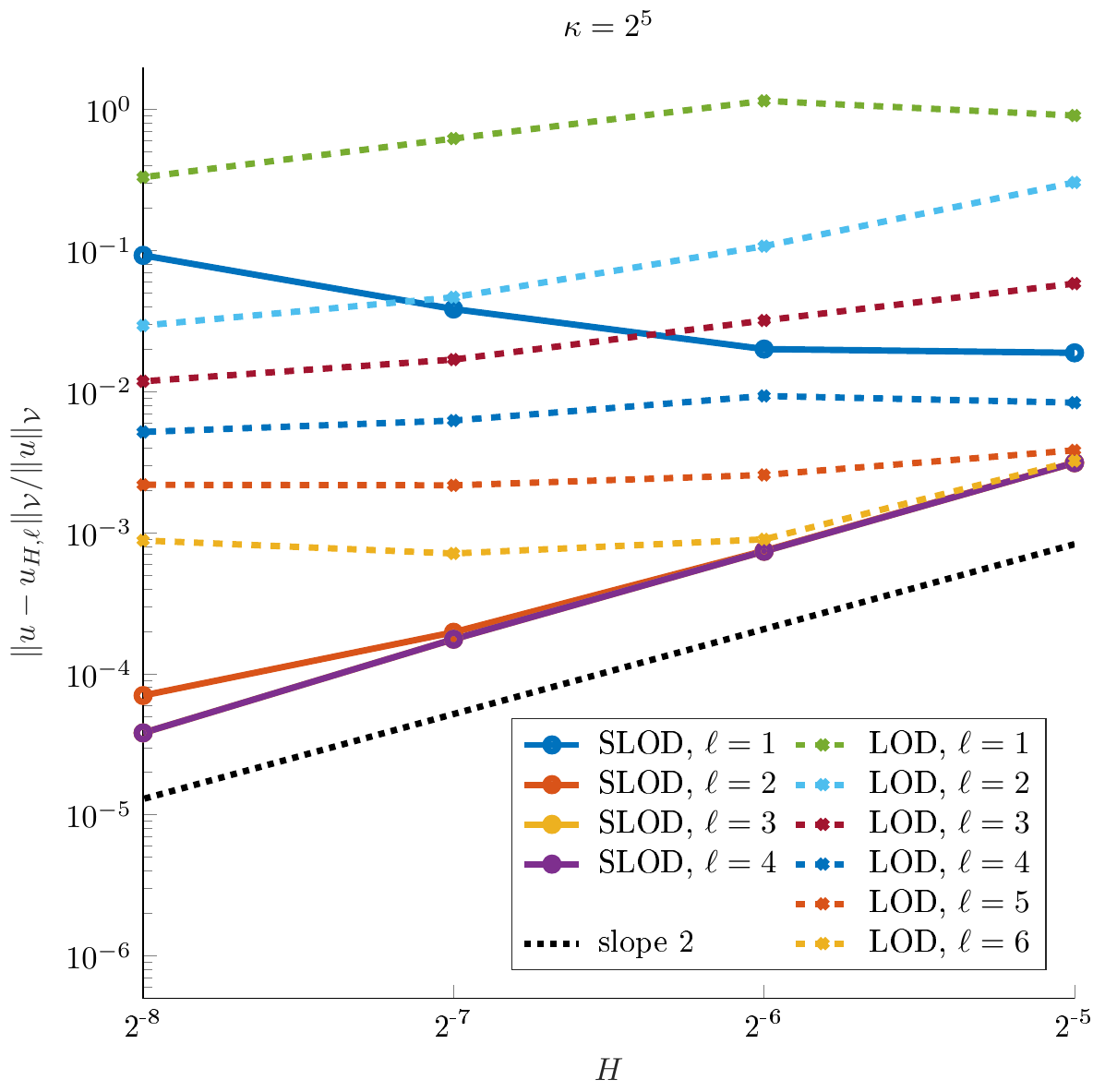}
	\includegraphics[width=.49\textwidth]{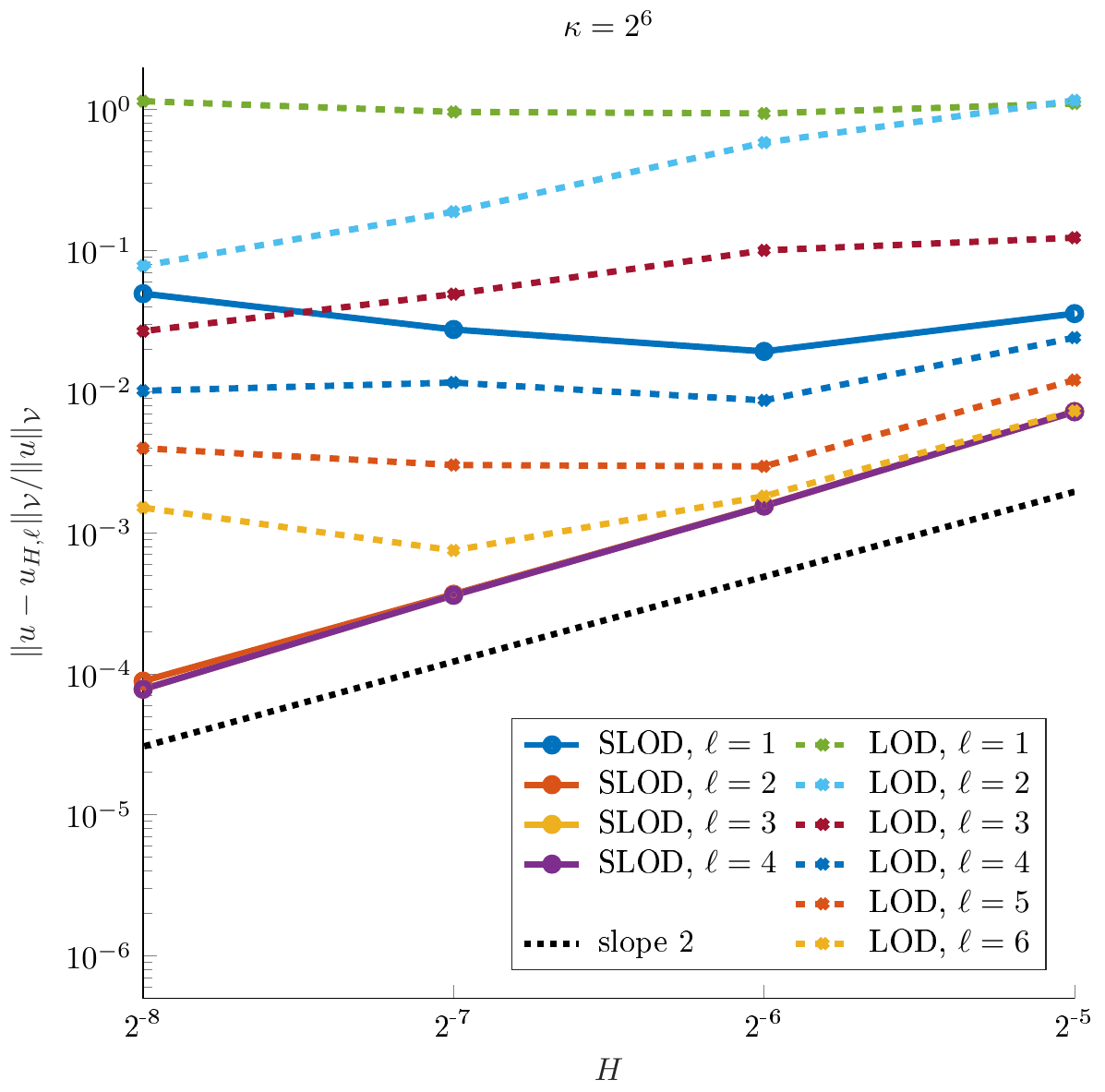}
	\caption{Convergence plots for the SLOD and LOD for several oversampling parameters $\ell$ for $\kappa = 2^5$ (left) and $\kappa = 2^6$ (right). }
	\label{fig:conv}
\end{figure}

Here, we consider the right-hand side $f(x,y) = \sin(\pi x)\cos(\pi y)$. Figure  \ref{fig:conv} depicts the relative $\V$-norm errors for several $\ell$ in dependence of $H$ for $\kappa =  2^5, 2^6$ in a double-logarithmic plot. As reference, we indicate a line of slope 2 which is the expected convergence rate in $H$ for right-hand sides $f \in H^1(\Omega)$; see Theorem \ref{th:conv}.

In Figure \ref{fig:conv}, one clearly observes convergence of order 2 in $H$ for the SLOD.  Note that for the SLOD, the lines corresponding to the oversampling parameters $\ell = 2,3$ can hardly be distinguished. As, for the LOD, the localization error decays much slower than for the SLOD (see Figure \ref{fig:loc}), the localization error dominates the overall error. Therefore, for the LOD, one rather observes the decay of the localization error than the desired convergence in $H$. 

\subsection{High-contrast heterogeneous media}
\begin{figure}[h]
	\includegraphics[width=.49\textwidth]{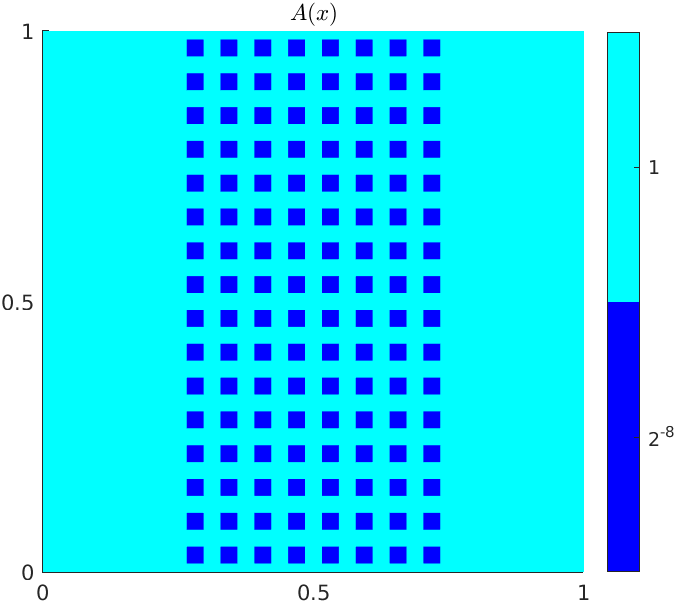}\hfill
	\includegraphics[width=.49\textwidth]{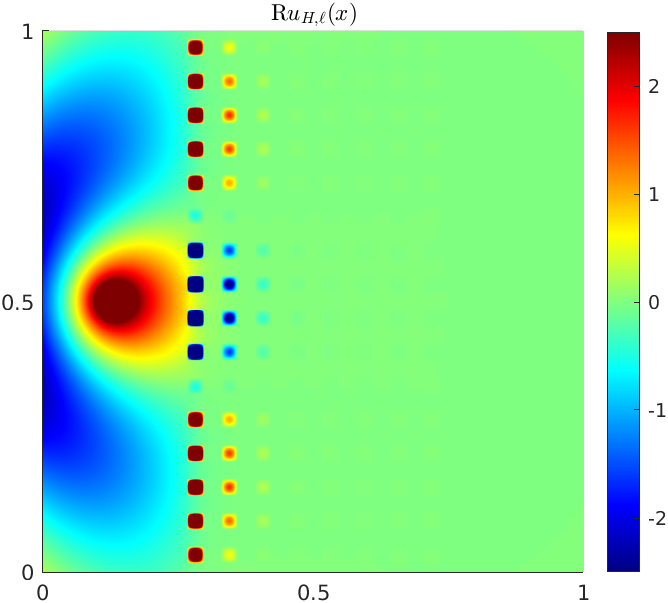}
	\caption{Heterogeneous coefficient $A$ (left) and real part of the corresponding SLOD solution (right) for $\ell = 2$, $H= 2^{-6}$, and $\kappa = 9$.}
	\label{fig:scatterer}
\end{figure}

Here, we demonstrate the application of the SLOD to the heterogeneous Helmholtz equation $-\nabla\cdot  (A\nabla u) -\kappa^2 u = f$ with homogeneous impedance boundary conditions. For some parameter $0<\varepsilon\ll1$, the coefficient $A$ takes the value $\varepsilon^2$ inside some periodically aligned inclusions of size $\varepsilon/2$ and the value $1$ elsewhere;  see Figure \ref{fig:scatterer} (left) for a depiction of $A$. For the choice $\kappa = 9$, a special interplay between the wavenumber and the periodic structure of the inclusions yields a negative-valued effective wavenumber in homogenization theory which triggers an exponential decay of the modulus of the Helmholtz solution in the bulk domain. This physically interesting effect is caused by Mie resonances in the small inclusions; see \cite{PeV20}. As right-hand side, we use an approximate  point source located at $z = (0.125,0.5)^T$ that vanishes outside a circle of radius $0.05$, i.e., 
\begin{align}\label{def:pointsource}
	f(x,y) = \begin{cases}
	10^4 \cdot  \exp\Bigg(\frac{-1}{1-\frac{(x - z_1)^2 + (y - z_2)^2}{0.05^2}}\Bigg),& (x - z_1)^2 + (y - z_2)^2 < 0.05^2 \\
	0,& \text{else}
	\end{cases}.
\end{align}
Figure \ref{fig:scatterer} (right) depicts the real part of the SLOD solution for $H = 2^{-6}$ and  $\ell = 2$. Note that, for the sake of illustration, the color map is truncated to the interval $[-2.5,2.5]$. The SLOD solution has a relative error of $3.3\cdot 10^{-3}$  with respect to the weighted norm  $$\|\cdot \|_{\V,A}^2\coloneqq \tnormf{A^{1/2}\nabla\cdot}{\Omega}^2 + \kappa^2\tnormf{\cdot}{\Omega}^2.$$
For reaching a similar accuracy on the same coarse mesh, the LOD needs an oversampling parameter of $\ell = 5$ which is a significant difference to $\ell = 2$ for the SLOD.

\subsection{Perfectly matched layers (PML)}\label{sec:Perfectly matched layers}
For this numerical experiment, we again consider the point source \eqref{def:pointsource}, but this time with  $z = (0.5,0.5)^T$. Our implementation of the perfectly matched layer is similar to \cite{PML1}, but adapted to the multi-scale setting.  We consider the fixed coarse Cartesian mesh $\TH$ with $H = 2^{-7}$ and divide the domain $\Omega$ into an inner part $\Omega_\mathrm{F}\coloneqq (4H,1-4H)^2$, the physical domain, and the absorbing layer $\Omega_\mathrm{A}\coloneqq \Omega \backslash\Omega_\mathrm{F}$. For this configuration, the absorbing layer has a width of $4H$. The (unbounded) absorbing function in $x$-direction is given as
\begin{align*}
	\rho_x(x) = \begin{cases}
		\frac{\mathrm{i}}{\kappa}\left(\frac{1}{-x} + \frac{1}{4H}\right), & 0 < x \leq 4H \\
		\frac{\mathrm{i}}{\kappa}\left(\frac{1}{1-x} - \frac{1}{4H}\right), & 1-4H \leq x < 1
	\end{cases}.
\end{align*}
In the $y$-direction the absorbing functions are chosen accordingly. As usual for PML, we use homogeneous Dirichlet boundary conditions on $\partial\Omega$. The full PML Helmholtz system may be found in \cite[Section 3]{PML1}. We apply the SLOD to this PML formulation and truncate the solution to the physical domain $\Omega_\mathrm{F}$. 

\begin{figure}[h]
	\includegraphics[width=.49\textwidth]{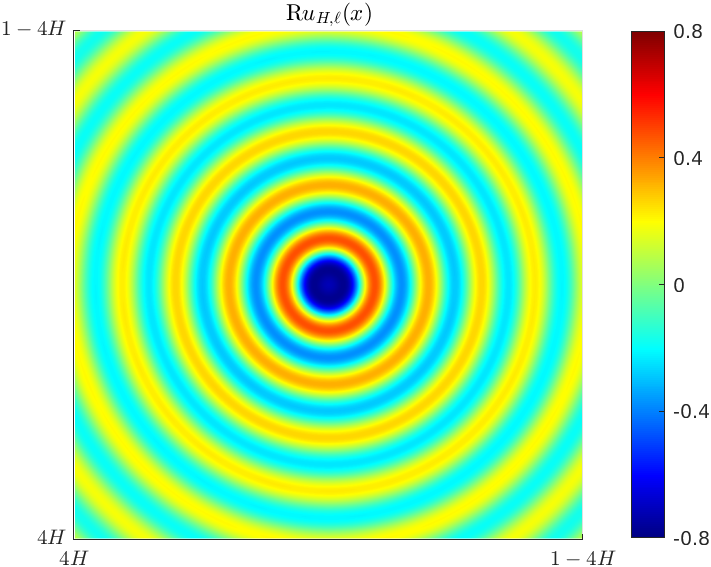}\hfill
	\includegraphics[width=.49\textwidth]{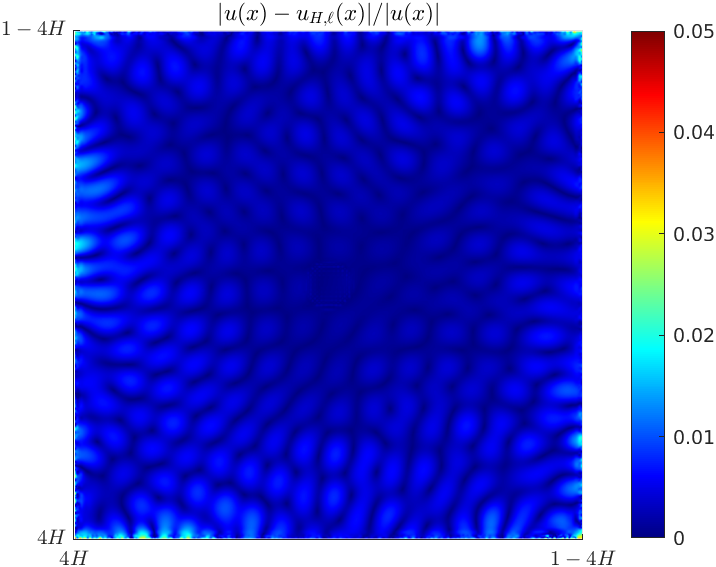}
	\caption{Real part of SLOD solution with PML (left) and relative error (right) for $\ell = 2$, $H = 2^{-7}$, and $\kappa = 2^6$.}
	\label{fig:PML}
\end{figure}

We choose the  parameters $\ell = 2$ and $\kappa = 2^6$. Figure \ref{fig:PML} shows the real part of the SLOD solution with PML (left) and the relative error computed against the PML reference solution (right). The relative error with respect to  $\|\cdot\|_{\V_{\Omega_\mathrm{F}}}$ is $6.5\cdot 10^{-3}$.

\section{Conclusion}
In this paper, we introduced a novel multi-scale method for high-frequency Helmholtz problems. It is conceptually similar to the LOD but utilizes a substantially improved localization strategy. The resulting relaxed oversampling condition has a significant impact on the computational costs; in the offline phase as well as in the online phase, significant savings can be achieved.

Under a stability assumption on the method's basis, a rigorous wavenumber-explicit stability and error analysis of the proposed method was performed. Additionally, we proposed an a-posteriori error control strategy for the localization error, which uses easily computable (quasi-)local error indicators.  It controls the localization error by adaptively increasing the patch-size of patches with a large associated error indicator.

A sequence of numerical experiments demonstrated the effectiveness of the proposed multi-scale method. In contrast to the LOD, it yields faithful  numerical approximations already for the relatively small oversampling parameters $\ell= 2,3$. In practice, the dependence of $\ell$ on the wavenumber is hardly noticeable.  We demonstrated that the proposed method can handle numerically challenging high-contrast heterogeneous Helmholtz problems. Furthermore, it was shown that the method is easily combined with perfectly matched layers (PML). We highlight that our approach may be transferred to other related problems, such as elastic wave propagation \cite{bg2016} or Maxwell's equations \cite{Gallistl18,HeP20} that have already been studied in the classical LOD setting.

\bibliographystyle{alpha}
\bibliography{bib}
\end{document}